\documentclass[12pt,english]{article}
\makeatletter

\usepackage{ifpdf}

\newif\ifpdf
\ifx\pdfoutput\undefined
  \pdffalse
\else
  \pdfoutput=1
  \pdftrue
\fi

\RequirePackage{xspace} %
\RequirePackage{subfigure} %

\ifpdf
  \RequirePackage[ pdftex, plainpages = false, pdfpagelabels,
                 pdfpagelayout = useoutlines,
                 bookmarks,
                 breaklinks = true,
                 linktocpage,
                 colorlinks = true,
                 linkcolor = blue,
                 urlcolor  = blue,
                 citecolor = blue,
                 anchorcolor = blue,
                 hyperindex = true,
                 hyperfigures
                 ]{hyperref}
\usepackage{makeidx} 
\usepackage{setspace} 
\usepackage{rotating} 
\usepackage{ecltree}
\usepackage{epic}
\usepackage{supertabular}  
\usepackage{color}
\usepackage{exscale}
\usepackage{fontenc}
\usepackage{ifthen}
\usepackage{latexsym}
\usepackage{makeidx}
\usepackage{syntonly}
\usepackage{inputenc}
\usepackage{graphicx}
\usepackage{setspace}
\usepackage{caption2}
\usepackage[english]{babel}
\usepackage[square, comma,numbers,sort&compress]{natbib}
\usepackage{hypernat}
\usepackage{boxedminipage}
\usepackage{framed}
\usepackage{longtable}
\usepackage[all]{hypcap}
\usepackage{amsmath}
\usepackage{amssymb}
\usepackage{stmaryrd}
\usepackage{setspace}
\usepackage[bottom]{footmisc}

\newcommand{\note}[1]{\marginpar[left]{\singlespace \tiny #1}}

\renewcommand{\sectionmark}[1]%
      {\markright{\thesection\ #1}} 

\renewcommand{\note}[1]{}
\renewcommand\@biblabel[1]{$[{#1}]$}

\begin{document}

\title{Tensor Calculus\\
\vspace{6cm}}

\author{Taha Sochi\thanks{Department of Physics \& Astronomy, University College London, Gower
Street, London, WC1E 6BT. Email: t.sochi@ucl.ac.uk.}\vspace{8cm}}

\maketitle
\pagebreak{}

\phantomsection \addcontentsline{toc}{section}{Preface}

\section*{Preface}

These notes are the second part of the tensor calculus documents which
started with the previous set of introductory notes \cite{SochiTensorIntro2016}.
In the present text, we continue the discussion of selected topics
of the subject at a higher level expanding, when necessary, some topics
and developing further concepts and techniques. The purpose of the
present text is to solidify, generalize, fill the gaps and make more
rigorous what have been presented in the previous set of notes and
to prepare the ground for the next set of notes. Unlike the previous
notes which are largely based on a Cartesian approach, the present
notes are essentially based on assuming an underlying general curvilinear
coordinate system. We also provide a small sample of proofs to familiarize
the reader with the tensor techniques inline with the tutorial nature
of the present text; however, due to the limited objectives of the
present text we do not provide comprehensive proofs and complete theoretical
foundations for the provided materials.

We generally follow the same conventions and notations used in the
previous set of notes with the following amendments:

$\bullet$ We use capital Gamma, $\Gamma_{jk}^{i}$, for the Christoffel
symbols of the second kind which is more elegant and readable than
the curly bracket notation $\left\{ _{jk}^{i}\right\} $ that we used
in the previous notes insisting that, despite the suggestive appearance
of the Gamma notation, the Christoffel symbols are not tensors in
general.

$\bullet$ Due to the restriction of using real (non-complex) quantities,
as stated in the previous notes, all arguments of real-valued functions,
like square roots and logarithmic functions, are assumed to be non-negative
by taking the absolute value, if necessary, without using the absolute
value symbol, as done by some authors. This is to simplify the notation
and avoid confusion with the determinant notation.

$\bullet$ We generalize the partial derivative notation so that $\partial_{i}$
can symbolize the partial derivative with respect to the $u^{i}$
coordinate of general curvilinear systems and not just for Cartesian
coordinates which are usually denoted by $x^{i}$. The type of coordinates,
being Cartesian or general or otherwise, will be determined by the
context which should be obvious in all cases.

$\bullet$ The summation symbol (i.e. $\sum$) is used in most cases
when a summation is needed but the summation convention conditions
do not apply or there is an ambiguity about it, e.g. when an index
is repeated more than twice or a twice-repeated index is in an upper
or lower state in both positions or a summation index is not repeated
visually because it is part of a squared symbol.

$\bullet$ ``Tensor'' and ``Matrix'' are not the same; however
for ease of expression they are used sometimes interchangeably and
hence some tensors may be referred to as matrices meaning the matrix
representing the tensor.

$\bullet$ In the present text, all coordinate transformations are
assumed to be continuous, single valued and invertible.

\pagebreak{}

\phantomsection \addcontentsline{toc}{section}{Contents}

\tableofcontents{}

\pagebreak{}

\section{Coordinate Systems, Spaces and Transformations\label{secCoordinateSystems}}

$\bullet$ The focus of this section is coordinate systems, their
types and transformations as well as some general properties of spaces
which are needed for the development of the concepts and techniques
of tensor calculus in the present and forthcoming notes.

\subsection{Coordinate Systems}

$\bullet$ In simple terms, a coordinate system is a mathematical
device, essentially of geometric nature, used by an observer to identify
the location of points and objects and describe events in generalized
space which may include space-time.

$\bullet$ The coordinates of a system can have the same or different
physical dimensions. An example of the first is the Cartesian system
where all the coordinates have the dimension of length, while examples
of the second include the cylindrical and spherical systems where
some coordinates have the dimension of length while others are dimensionless.

$\bullet$ Generally, the physical dimensions of the components and
basis vectors of the covariant and contravariant forms of a tensor
are different.

\subsection{Spaces}

$\bullet$ A Riemannian space is a manifold characterized by the existing
of a symmetric rank-2 tensor called the metric tensor. The components
of this tensor, which can be in covariant ($g_{ij}$) or contravariant
($g^{ij}$) forms, are in general continuous variable functions of
coordinates, i.e. $g_{ij}=g_{ij}(u^{1},u^{2},\ldots,u^{n})$ and $g^{ij}=g^{ij}(u^{1},u^{2},\ldots,u^{n})$
where $u^{i}$ symbolize general coordinates. This tensor facilitates,
among other things, the generalization of lengths and distances in
general coordinates where the length of an element of arc, $ds$,
is defined by:
\begin{equation}
\left(ds\right)^{2}=g_{ij}du^{i}du^{j}
\end{equation}
In the special case of a Euclidean space coordinated by a rectangular
system, the metric becomes the identity tensor, that is:
\begin{equation}
g_{ij}=g^{ij}=g_{j}^{i}=\delta_{ij}=\delta^{ij}=\delta_{j}^{i}
\end{equation}

$\bullet$ The metric of a Riemannian space may be called the Riemannian
metric. Similarly, the geometry of the space may be described as a
Riemannian geometry.

$\bullet$ All spaces dealt with in the present notes are Riemannian
with well-defined metrics.

$\bullet$ A manifold or space is dubbed ``flat'' when it is possible
to find a coordinate system for the space with a diagonal metric tensor
whose all diagonal elements are $\pm1$; the space is called ``curved''
otherwise. Examples of flat space are the 3D Euclidean space coordinated
by a rectangular Cartesian system whose metric tensor is diagonal
with all the diagonal elements being $+1$, and the 4D Minkowski space-time
whose metric is diagonal with elements of $\pm1$. Examples of curved
space is the 4D space-time of general relativity in the presence of
matter and energy.

$\bullet$ When all the diagonal elements of the metric tensor of
a flat space are $+1$, the space and the coordinate system may be
described as homogeneous.

$\bullet$ An $n$D manifold is Euclidean \textit{iff} $R_{ijkl}=0$
where $R_{ijkl}$ is the Riemann tensor (see $\S$ \ref{subRiemannTensor});
otherwise the manifold is curved to which the general Riemannian geometry
applies.

$\bullet$ A ``field'' is a function of the position vector over
a region of space. Scalars, vectors and tensors may be defined on
a single point of the space or over an extended region of the space;
in the latter case we have scalar fields, vector fields and tensor
fields, e.g. temperature field, velocity field and stress field respectively.

$\bullet$ In metric spaces, the physical quantities are independent
of the form of description, being covariant or contravariant, as the
metric tensor facilitates the transformation between the different
forms; hence making the description objective.

\subsection{Transformations}

$\bullet$ In general terms, a transformation from an $n$D space
to another $n$D space is a correlation that maps a point from the
first space (original) to a point in the second space (transformed)
where each point in the original and transformed spaces is identified
by $n$ independent variables or coordinates. To distinguish between
the two sets of coordinates in the two spaces, the coordinates of
the points in the transformed space may be notated with barred symbols,
e.g. ($\bar{u}^{1},\bar{u}^{2},\ldots,\bar{u}^{n}$) or ($\bar{u}_{1},\bar{u}_{2},\ldots,\bar{u}_{n}$)
where the superscripts and subscripts are indices, while the coordinates
of the points in the original space are notated with unbarred similar
symbols, e.g. ($u^{1},u^{2},\ldots,u^{n}$) or ($u_{1},u_{2},\ldots,u_{n}$).
Under certain conditions, which will be clarified later, such a transformation
is unique and hence an inverse transformation from the transformed
space to the original space is also defined. Mathematically, each
one of the direct and inverse transformations can be regarded as a
correlation expressed by a set of equations in which each coordinate
in one space is considered as a function of the coordinates in the
other space. Hence the transformations between the two sets of coordinates
in the two spaces can by expressed mathematically by the following
two sets of independent relations:
\begin{equation}
\bar{u}^{i}=\bar{u}^{i}(u^{1},u^{2},\ldots,u^{n})\,\,\,\,\,\,\,\,\,\,\,\,\,\,\,\,\,\,\,\&\,\,\,\,\,\,\,\,\,\,\,\,\,\,\,\,\,\,u^{i}=u^{i}(\bar{u}^{1},\bar{u}^{2},\ldots,\bar{u}^{n})\label{eqTransformationEquations}
\end{equation}
where $i=1,2,\ldots,n$ with $n$ being the space dimension. The independence
of the above relations is guaranteed \textit{iff} the Jacobian of
the transformation does not vanish on any point in the space (see
about Jacobian the forthcoming points). An alternative to viewing
the transformation as a mapping between two different spaces is to
view it as a correlation of the same point in the same space but observed
from two different coordinate frames of reference which are subject
to a similar transformation. The following points will be largely
based on the latter view.

$\bullet$ As far as the notation is concerned, there is no fundamental
difference between the barred and unbarred systems and hence the notation
can be interchanged.

$\bullet$ An injective transformation maps any two distinct points
of the original space, $\mathbf{r}_{1}$ and $\mathbf{r}_{2}$, onto
two distinct points of the transformed space, $\mathbf{\bar{r}}_{1}$
and $\bar{\mathbf{r}}_{2}$. The image of an injective transformation,
$\mathbf{\bar{r}}$, is regarded as coordinates for the point, and
the collection of all such coordinates of the space points may be
considered as a representation of a coordinate system for the space.
If the mapping from an original rectangular system is linear, the
coordinate system obtained from such a transformation is called ``affine''.
Coordinate systems which are not affine are described as ``curvilinear''
such as cylindrical and spherical systems.

$\bullet$ The following $n\times n$ matrix of $n^{2}$ partial derivatives
of the barred coordinates with respect to the unbarred coordinates
is called the ``Jacobian matrix'' of the transformation between
the barred and unbarred systems:
\begin{equation}
\mathbf{J}=\left[\begin{array}{cccc}
\frac{\partial\bar{u}^{1}}{\partial u^{1}} & \frac{\partial\bar{u}^{1}}{\partial u^{2}} & \cdots & \frac{\partial\bar{u}^{1}}{\partial u^{n}}\\
\frac{\partial\bar{u}^{2}}{\partial u^{1}} & \frac{\partial\bar{u}^{2}}{\partial u^{2}} & \cdots & \frac{\partial\bar{u}^{2}}{\partial u^{n}}\\
\vdots & \vdots & \ddots & \vdots\\
\frac{\partial\bar{u}^{n}}{\partial u^{1}} & \frac{\partial\bar{u}^{n}}{\partial u^{2}} & \cdots & \frac{\partial\bar{u}^{n}}{\partial u^{n}}
\end{array}\right]
\end{equation}
while its determinant:
\begin{equation}
J=\mathrm{det}(\mathbf{J})
\end{equation}
is called the ``Jacobian'' of the transformation.

$\bullet$ Barred and unbarred in the definition of Jacobian should
be understood in a general sense not just as two labels since the
Jacobian is not restricted to transformations between two systems
of the same type but labeled as barred and unbarred. In fact the two
coordinate systems can be fundamentally different in nature and not
of the same type such as Cartesian and general curvilinear. The Jacobian
matrix and determinant represent any transformation by the above partial
derivative matrix system between two coordinates defined by two different
sets of coordinate variables not necessarily as barred and unbarred.
The objective of defining the Jacobian as between barred and unbarred
systems is generality and clarity.

$\bullet$ The transformation from the unbarred coordinate system
to the barred coordinate system is bijective\footnote{Bijective transformation means injective (one-to-one) and surjective
(onto) mapping.} \textit{iff} $J\ne0$ on any point in the transformed region of the
space. In this case the inverse transformation from the barred to
the unbarred system is also defined and bijective and is represented
by the inverse of the Jacobian matrix:\footnote{Notationally, there is no fundamental difference between the barred
and unbarred systems and hence the labeling is rather arbitrary and
can be interchanged. Therefore, the Jacobian may be notated as barred
over unbarred or the other way around. Yes, in a specific context
when one of these is labeled as the Jacobian, the other one should
be labeled as the inverse Jacobian to distinguish between the two
opposite Jacobians and their corresponding transformations.}
\begin{equation}
\mathbf{\bar{J}}=\mathbf{J}^{-1}
\end{equation}
Consequently, the Jacobian of the inverse transformation, being the
determinant of the inverse Jacobian matrix, is the reciprocal of the
Jacobian of the original transformation:
\begin{equation}
\bar{J}=\frac{1}{J}
\end{equation}

$\bullet$ A coordinate transformation is admissible \textit{iff}
the transformation is bijective with non-vanishing Jacobian and the
transformation function is of class $C^{2}$.\footnote{$C^{n}$ continuity condition means that the function and all its
first $n$ partial derivatives do exist and are continuous in their
domain of definition. Also, some authors impose a weaker condition
of being of class $C^{1}$.}

$\bullet$ ``Affine tensors'' are tensors that correspond to admissible
linear coordinate transformations from an original rectangular system
of coordinates.

$\bullet$ Coordinate transformations are described as ``proper''
when they preserve the handedness (right- or left-handed) of the coordinate
system and ``improper'' when they reverse the handedness. Improper
transformations involve an odd number of coordinate axes inversions
through the origin.

$\bullet$ Inversion of axes may be called improper rotation while
ordinary rotation is described as proper rotation.

$\bullet$ Transformations of coordinates can be active, when they
change the state of the observed object such as rotating the object
in the space, or passive when they are based on keeping the state
of the object and changing the state of the coordinate system which
the object is observed from. In brief, the subject of an active transformation
is the object while the subject of a passive transformation is the
coordinate system.

$\bullet$ An object that does not change by admissible coordinate
transformations is described as ``invariant'' such as the value
of a true scalar and the length of a vector.

$\bullet$ As there are essentially two different types of basis vectors,
namely tangent vectors of covariant nature and gradient vectors of
contravariant nature, there are two main types of non-scalar tensors:
contravariant and covariant tensors which are based on the type of
the employed basis vectors. Tensors of mixed type employ in their
definition mixed basis vectors of the opposite type to the corresponding
indices of their components. As indicated earlier, the transformation
between these different types is facilitated by the metric tensor.

$\bullet$ A product or composition of coordinate transformations
is a succession of transformations where the output of one transformation
is taken as the input to the next transformation. In such cases, the
Jacobian of the product is the product of the Jacobians of the individual
transformations of which the product is made.

$\bullet$ The collection of all admissible coordinate transformations
with non-vanishing Jacobian form a group, that is they satisfy the
properties of closure, associativity, identity and inverse. Hence,
any convenient admissible coordinate system can be chosen as the point
of entry since other systems can be reached, if needed, through the
set of admissible transformations. This is the cornerstone of building
covariant physical theories which are independent of the subjective
choice of coordinate systems and reference frames.

$\bullet$ Transformation of coordinates is not a commutative operation.

\subsection{Coordinate Surfaces and Curves}

$\bullet$ The surfaces of constant coordinates at a certain point
of the space meet to form curves (i.e. curves of intersection of these
surfaces in pairs). The coordinate curves are these curves of mutual
intersection of the surfaces of constant coordinates of the curvilinear
system.

$\bullet$ The above transformation equations (Eq. \ref{eqTransformationEquations})
are used to define the set of surfaces of constant coordinates and
coordinate curves of mutual intersection of these surfaces. These
coordinate surfaces and curves play a crucial role in the formulation
and development of this subject.

$\bullet$ The coordinate axes of a coordinate system can be rectilinear,
and hence the coordinate curves are straight lines and the surfaces
of constant coordinates are planes, as in the case of rectangular
Cartesian systems, or curvilinear, and hence the coordinate curves
are generalized curved paths and the surfaces of constant coordinates
are generalized curved surfaces, as in the case of cylindrical and
spherical systems.

$\bullet$ In curvilinear coordinate systems, some or all of the coordinate
surfaces are not planes and some or all of the coordinate lines are
not straight lines.

$\bullet$ Orthogonal coordinate systems are those for which the vectors
tangent to the coordinate curves, as well as the vectors normal to
the surfaces of constant coordinates, are mutually perpendicular at
all points of the space. Consequently, in orthogonal coordinates,
the coordinate surfaces are mutually perpendicular and the coordinate
lines are also perpendicular at the point of intersection.

$\bullet$ In orthogonal coordinate systems, the corresponding covariant
and contravariant basis vectors at any given point in the space are
in the same direction, i.e. the tangent vector to a particular coordinate
curve $u^{i}$ at a certain point and the gradient vector normal to
the surface of constant $u^{i}$ at the same point have the same direction
although they may be of different length.

$\bullet$ A necessary and sufficient condition for a coordinate system
to be orthogonal is that its metric tensor is diagonal.

$\bullet$ An admissible coordinate transformation from a Cartesian
system defines another Cartesian system if the transformation is linear,
and defines a curvilinear system if the transformation is nonlinear.

\subsection{Scale Factors}

$\bullet$ Scale factors (usually symbolized with $h_{1},h_{2},\ldots h_{n}$)
of a coordinate system are those factors which are required to multiply
the coordinate differentials to obtain distances traversed during
a change in the coordinate of that magnitude, e.g. $\rho$ in the
plane polar coordinate system which multiplies the differential of
the polar angle $d\phi$ to obtain the distance $L$ traversed by
a change of magnitude $d\phi$ in the polar angle which is $L=\rho\,d\phi$.
They are also used to normalize the basis vectors (refer to the previous
notes \cite{SochiTensorIntro2016} and forthcoming notes).

$\bullet$ The scale factors for the Cartesian, cylindrical and spherical
coordinate systems in 3D spaces are given in Table \ref{tabScaleFactors}.

\begin{table}
\centering{}\caption{The scale factors for the three most commonly used orthogonal coordinate
systems in 3D spaces. The squares of these entries and the reciprocals
of these squares give the diagonal elements of the covariant and contravariant
metric tensors $g_{ij}$ and $g^{ij}$ respectively of these systems.\vspace{0.2cm}\label{tabScaleFactors}}
\begin{tabular*}{16cm}{@{\extracolsep{\fill}}|c|c|c|c|}
\hline
 & Cartesian ($x,y,z$)~~~~~~ & Cylindrical ($\rho,\phi,z$)~~~~~~ & Spherical ($r,\theta,\phi$)~~~~~~\tabularnewline
\hline
$h_{1}$ & 1 & 1 & 1\tabularnewline
\hline
$h_{2}$ & 1 & $\rho$ & $r$\tabularnewline
\hline
$h_{3}$ & 1 & 1 & $r\,\sin\theta$\tabularnewline
\hline
\end{tabular*}
\end{table}

$\bullet$ The scale factors are also used in the expressions for
the differential elements of arc, surface and volume in general orthogonal
coordinates, as described in $\S$ \ref{subMetricTensor}.

\subsection{Basis Vectors in General Curvilinear Systems\label{subBasisVectors}}

$\bullet$ The vectors providing the basis set for a coordinate system,
which are not necessarily of unit length or mutually orthogonal, are
of covariant type when they are tangent to the coordinate curves,
and of contravariant type when they are perpendicular to the local
surfaces of constant coordinates. Formally, the covariant and contravariant
basis vectors are defined respectively by:
\begin{equation}
\mathbf{E}_{i}=\frac{\partial\mathbf{r}}{\partial u^{i}}\,\,\,\,\,\,\,\,\,\,\,\,\,\,\,\,\,\,\,\,\&\,\,\,\,\,\,\,\,\,\,\,\,\,\,\,\,\,\,\,\,\mathbf{E}^{i}=\nabla u^{i}
\end{equation}
where $\mathbf{r}$ is the position vector in Cartesian coordinates
($x^{1},x^{2},\ldots$), and $u^{i}$ are generalized curvilinear
coordinates.

$\bullet$ In general curvilinear coordinate systems, the covariant
and contravariant basis sets, $\mathbf{E}_{i}$ and $\mathbf{E}^{i}$,
are functions of coordinates, i.e.
\begin{equation}
\mathbf{E}_{i}=\mathbf{E}_{i}\left(u^{1},\ldots,u^{n}\right)\,\,\,\,\,\,\,\,\,\,\,\,\,\,\,\,\&\,\,\,\,\,\,\,\,\,\,\,\,\,\,\mathbf{E}^{i}=\mathbf{E}^{i}\left(u^{1},\ldots,u^{n}\right)
\end{equation}

$\bullet$ Like other vectors, the covariant and contravariant basis
vectors are related to each other through the metric tensor, that
is:
\begin{equation}
\mathbf{E}_{i}=g_{ij}\mathbf{E}^{j}\,\,\,\,\,\,\,\,\,\,\,\,\,\,\,\,\,\,\,\,\,\,\,\&\,\,\,\,\,\,\,\,\,\,\,\,\,\,\,\,\,\,\,\,\mathbf{E}^{i}=g^{ij}\mathbf{E}_{j}
\end{equation}

$\bullet$ The covariant and contravariant basis vectors are reciprocal
basis systems, and hence in a 3D space with a right-handed coordinate
system ($u_{1},u_{2},u_{3}$) they are linked by the following relations:
\begin{equation}
\mathbf{E}^{1}=\frac{\mathbf{E}_{2}\times\mathbf{E}_{3}}{\mathbf{E}_{1}\cdot\left(\mathbf{E}_{2}\times\mathbf{E}_{3}\right)},\,\,\,\,\,\,\,\,\,\,\,\,\,\,\mathbf{E}^{2}=\frac{\mathbf{E}_{3}\times\mathbf{E}_{1}}{\mathbf{E}_{1}\cdot\left(\mathbf{E}_{2}\times\mathbf{E}_{3}\right)},\,\,\,\,\,\,\,\,\,\,\,\,\,\,\mathbf{E}^{3}=\frac{\mathbf{E}_{1}\times\mathbf{E}_{2}}{\mathbf{E}_{1}\cdot\left(\mathbf{E}_{2}\times\mathbf{E}_{3}\right)}
\end{equation}

\begin{equation}
\mathbf{E}_{1}=\frac{\mathbf{E}^{2}\times\mathbf{E}^{3}}{\mathbf{E}^{1}\cdot\left(\mathbf{E}^{2}\times\mathbf{E}^{3}\right)},\,\,\,\,\,\,\,\,\,\,\,\,\,\,\mathbf{E}_{2}=\frac{\mathbf{E}^{3}\times\mathbf{E}^{1}}{\mathbf{E}^{1}\cdot\left(\mathbf{E}^{2}\times\mathbf{E}^{3}\right)},\,\,\,\,\,\,\,\,\,\,\,\,\,\,\mathbf{E}_{3}=\frac{\mathbf{E}^{1}\times\mathbf{E}^{2}}{\mathbf{E}^{1}\cdot\left(\mathbf{E}^{2}\times\mathbf{E}^{3}\right)}
\end{equation}

$\bullet$ The relations in the last point may be expressed in a more
compact form as follow:
\begin{equation}
\mathbf{E}^{i}=\frac{\mathbf{E}_{j}\times\mathbf{E}_{k}}{\mathbf{E}_{i}\cdot\left(\mathbf{E}_{j}\times\mathbf{E}_{k}\right)}\,\,\,\,\,\,\,\,\,\,\,\,\,\,\,\,\,\&\,\,\,\,\,\,\,\,\,\,\,\,\,\,\,\,\,\,\,\,\mathbf{E}_{i}=\frac{\mathbf{E}^{j}\times\mathbf{E}^{k}}{\mathbf{E}^{i}\cdot\left(\mathbf{E}^{j}\times\mathbf{E}^{k}\right)}
\end{equation}
where $i,j,k$ take respectively the values $1,2,3$ and the other
two cyclic permutations (i.e. $2,3,1$ and $3,1,2$).

$\bullet$ The magnitude of the scalar triple product $\mathbf{E}_{i}\cdot\left(\mathbf{E}_{j}\times\mathbf{E}_{k}\right)$
represents the volume of the parallelepiped formed by $\mathbf{E}_{i}$,
$\mathbf{E}_{j}$ and $\mathbf{E}_{k}$.

$\bullet$ The magnitudes of the basis vectors in general orthogonal
coordinates are given by:
\begin{equation}
\left|\mathbf{E}_{i}\right|=h_{i}\,\,\,\,\,\,\,\,\,\,\,\,\,\,\,\,\,\,\,\,\&\,\,\,\,\,\,\,\,\,\,\,\,\,\,\,\,\,\,\,\,\left|\mathbf{E}^{i}\right|=\frac{1}{h_{i}}
\end{equation}
where $h_{i}$ is the scale factor for the $i^{th}$ coordinate.

$\bullet$ The base vectors in the barred and unbarred general curvilinear
coordinate systems are related by the following transformation rules:
\begin{equation}
\begin{aligned}\mathbf{E}_{i} & =\frac{\partial\bar{u}^{j}}{\partial u^{i}}\bar{\mathbf{E}}_{j}\,\,\,\,\,\,\,\,\,\,\,\,\,\,\,\,\,\&\,\,\,\,\,\,\,\,\,\,\,\,\,\,\,\,\,\,\bar{\mathbf{E}}_{i}=\frac{\partial u^{j}}{\partial\bar{u}^{i}}\mathbf{E}_{j}\\
\mathbf{E}^{i} & =\frac{\partial u^{i}}{\partial\bar{u}^{j}}\bar{\mathbf{E}}^{j}\,\,\,\,\,\,\,\,\,\,\,\,\,\,\,\,\,\&\,\,\,\,\,\,\,\,\,\,\,\,\,\,\,\,\,\,\bar{\mathbf{E}}^{i}=\frac{\partial\bar{u}^{i}}{\partial u^{j}}\mathbf{E}^{j}
\end{aligned}
\end{equation}
where the indexed $u$ and $\bar{u}$ represent the coordinates in
the unbarred and barred systems respectively. The transformation rules
for the components can be straightforwardly concluded from the above
rules; for example for a vector $\mathbf{A}$ which can be represented
covariantly and contravariantly in the unbarred and barred systems
as:
\begin{equation}
\begin{aligned}\mathbf{A} & =\mathbf{E}^{i}A_{i}=\bar{\mathbf{E}}^{i}\bar{A}_{i}\\
\mathbf{A} & =\mathbf{E}_{i}A^{i}=\bar{\mathbf{E}}_{i}\bar{A}^{i}
\end{aligned}
\end{equation}
the transformation equations of its components between the two systems
are given respectively by:
\begin{equation}
\begin{aligned}A_{i} & =\frac{\partial\bar{u}^{j}}{\partial u^{i}}\bar{A}_{j}\,\,\,\,\,\,\,\,\,\,\,\,\,\,\,\,\,\&\,\,\,\,\,\,\,\,\,\,\,\,\,\,\,\,\,\,\bar{A}_{i}=\frac{\partial u^{j}}{\partial\bar{u}^{i}}A_{j}\\
A^{i} & =\frac{\partial u^{i}}{\partial\bar{u}^{j}}\bar{A}^{j}\,\,\,\,\,\,\,\,\,\,\,\,\,\,\,\,\,\&\,\,\,\,\,\,\,\,\,\,\,\,\,\,\,\,\,\,\bar{A}^{i}=\frac{\partial\bar{u}^{i}}{\partial u^{j}}A^{j}
\end{aligned}
\end{equation}
These transformation rules can be easily extended to higher rank tensors
of different variance types, as detailed in the introductory notes
\cite{SochiTensorIntro2016}.

$\bullet$ For a 3D manifold with a right-handed curvilinear coordinate
system, we have:
\begin{equation}
\mathbf{E}_{1}\cdot\left(\mathbf{E}_{2}\times\mathbf{E}_{3}\right)=\sqrt{g}\,\,\,\,\,\,\,\,\,\,\,\,\,\,\,\,\,\,\,\&\,\,\,\,\,\,\,\,\,\,\,\,\,\,\,\,\,\,\,\mathbf{E}^{1}\cdot\left(\mathbf{E}^{2}\times\mathbf{E}^{3}\right)=\frac{1}{\sqrt{g}}
\end{equation}
where $g$ is the determinant of the covariant metric tensor, i.e.
\begin{equation}
g=\mathrm{det}\left(g_{ij}\right)=\left|g_{ij}\right|
\end{equation}

$\bullet$ Because $\mathbf{E}_{i}\cdot\mathbf{E}_{j}=g_{ij}$ (Eq.
\ref{eqEg}) we have:
\begin{equation}
\mathbf{J}^{T}\mathbf{J}=\left[g_{ij}\right]\label{eqeqJTJgij}
\end{equation}
where $\mathbf{J}$ is the Jacobian matrix transforming between Cartesian
and generalized coordinates, the superscript $T$ represents matrix
transposition, $\left[g_{ij}\right]$ is the matrix representing the
covariant metric tensor and the product on the left is a matrix product
as defined in linear algebra which is equivalent to a dot product
in tensor algebra.

$\bullet$ Considering Eq. \ref{eqeqJTJgij}, the relation between
the determinant of the metric tensor and the Jacobian is given by:
\begin{equation}
g=J^{2}\label{eqgJsquared}
\end{equation}
where $J$ ($=\left|\frac{\partial x}{\partial u}\right|$ with $x$
for Cartesian and $u$ for generalized coordinates) is the Jacobian
of the transformation.

$\bullet$ As explained earlier, in orthogonal coordinate systems
the covariant and contravariant basis vectors, $\mathbf{E}_{i}$ and
$\mathbf{E}^{i}$, at any specific point of the space are in the same
direction, and hence the normalization of each one of these basis
sets, by dividing each basis vector by its magnitude, produces identical
orthonormal basis sets.\footnote{Consequently, there is no difference between the covariant and contravariant
components of tensors with respect to such contravariant and covariant
orthonormal basis sets.} This, however, is not true in general curvilinear coordinates where
each normalized basis set is different in general from the other.

$\bullet$ When the covariant basis vectors $\mathbf{E}_{i}$ are
mutually orthogonal at all points of the space, we have:

(A) the contravariant basis vectors $\mathbf{E}^{i}$ are mutually
orthogonal as well,

(B) the covariant and contravariant metric tensors, $g_{ij}$ and
$g^{ij}$, are diagonal with non-vanishing diagonal elements, i.e.
\begin{equation}
g_{ij}=g^{ij}=0\,\,\,\,\,\,\,\,\,\,\,\,\,\,\,\,\,\,\,\,\,\,\,\,\text{(\ensuremath{i\ne j})}
\end{equation}
\begin{equation}
g_{ii}\ne0\,\,\,\,\,\,\,\,\,\,\,\,\,\,\&\,\,\,\,\,\,\,\,\,\,\,\,\,g^{ii}\ne0\,\,\,\,\,\,\,\,\,\,\,\,\,\,\,\,\,\,\,\,\,\,\,\,\text{(no sum on \ensuremath{i})}
\end{equation}

(C) the diagonal elements of the covariant and contravariant metric
tensors are reciprocals, i.e.\footnote{The comment ``no summation'' may not be needed in this type of expressions
since both indices are of the same variance type in a generally non-Cartesian
system.}
\begin{equation}
g^{ii}=\frac{1}{g_{ii}}\,\,\,\,\,\,\,\,\,\,\,\,\,\,\,\,\,\,\,\text{(no summation)}
\end{equation}

(D) the magnitude of the contravariant and covariant basis vectors
are reciprocals, i.e.
\begin{equation}
\left|\mathbf{E}^{i}\right|=\frac{1}{\left|\mathbf{E}_{i}\right|}
\end{equation}

\subsection{Covariant, Contravariant and Physical Representations\label{subCovariantContravariantPhysical}}

$\bullet$ So far we are familiar with the covariant and contravariant
(including mixed) representations of tensors. There is still another
type of representation, that is the physical representation which
is the common one in the applications of tensor calculus such as fluid
and continuum mechanics.

$\bullet$ The covariant and contravariant basis vectors, as well
as the covariant and contravariant components of a vector, do not
in general have the same physical dimensions as indicated earlier;
moreover, the basis vectors may not have the same magnitude. This
motivates the introduction of a more standard form of vectors by using
physical components (which have the same dimensions) with normalized
basis vectors (which are dimensionless with unit magnitude) where
the metric tensor and the scale factors are employed to facilitate
this process. The normalization of the basis vectors is done by dividing
each vector by its magnitude. For example, the normalized covariant
basis vectors of a general coordinate system, $\hat{\mathbf{E}}_{i}$,
are given by:
\begin{equation}
\hat{\mathbf{E}}_{i}=\frac{\mathbf{E}_{i}}{\left|\mathbf{E}_{i}\right|}\,\,\,\,\,\,\,\,\,\,\,\,\,\,\,\,\,\,\,\,\,\,\,\,\,\,\,\,\,\,\,\,\text{(no sum on \ensuremath{i})}
\end{equation}
which for an orthogonal coordinate system becomes:
\begin{equation}
\hat{\mathbf{E}}_{i}=\frac{\mathbf{E}_{i}}{\sqrt{g_{ii}}}=\frac{\mathbf{E}_{i}}{h_{i}}\,\,\,\,\,\,\,\,\,\,\,\,\,\,\,\,\text{(no sum on \ensuremath{i})}
\end{equation}
where $g_{ii}$ is the $i^{th}$ diagonal element of the covariant
metric tensor and $h_{i}$ is the scale factor of the $i^{th}$ coordinate
as described previously. Consequently if the physical components of
a vector are notated with a hat, then for an orthogonal system we
have:
\begin{equation}
\mathbf{A}=A^{i}\mathbf{E}_{i}=\hat{A}^{i}\mathbf{\hat{E}}_{i}=\hat{A}^{i}\frac{\mathbf{E}_{i}}{\sqrt{g_{ii}}}\,\,\,\,\,\,\,\,\,\,\,\,\,\,\Longrightarrow\,\,\,\,\,\,\,\,\,\,\,\,\,\,\,\,\,\hat{A}^{i}=\sqrt{g_{ii}}A^{i}=h_{i}A^{i}\,\,\,\,\,\,\,\,\,\,\,\,\,\,\,\,\text{(no sum)}\label{eqPhysCompCova}
\end{equation}
Similarly for the contravariant basis vectors we have:
\begin{equation}
\mathbf{A}=A_{i}\mathbf{E}^{i}=\hat{A}_{i}\mathbf{\hat{E}}^{i}=\hat{A}_{i}\frac{\mathbf{E}^{i}}{\sqrt{g^{ii}}}\,\,\,\,\,\,\,\,\,\,\,\,\,\,\Longrightarrow\,\,\,\,\,\,\,\,\,\,\,\,\,\,\,\,\,\hat{A}_{i}=\sqrt{g^{ii}}A_{i}=\frac{A_{i}}{\sqrt{g_{ii}}}=\frac{A_{i}}{h_{i}}\,\,\,\,\,\,\,\,\,\,\,\,\,\,\,\,\text{(no sum)}
\end{equation}
where $g^{ii}$ is the $i^{th}$ diagonal element of the contravariant
metric tensor. These definitions and processes can be easily extended
to tensors of higher ranks.

$\bullet$ The physical components of higher rank tensors are similarly
defined as for rank-1 tensors by considering the basis vectors of
the coordinated space where similar simplifications apply to orthogonal
systems with mutually-perpendicular basis vectors. For example, for
a rank-2 tensor $\mathbf{A}$ with an orthogonal coordinate system,
the physical components can be represented by:
\begin{equation}
\begin{aligned}\hat{A}_{ij} & =\frac{A_{ij}}{h_{i}h_{j}} & \,\,\,\,\,\,\,\,\,\,\,\,\,\,\,\,\,\, & \text{(no sum on \ensuremath{i} or \ensuremath{j}, with basis \ensuremath{\hat{\mathbf{E}}^{i}\mathbf{\hat{E}}^{j}})}\\
\hat{A}^{ij} & =h_{i}h_{j}A^{ij} &  & \text{(no sum on \ensuremath{i} or \ensuremath{j}, with basis \ensuremath{\mathbf{\hat{E}}_{i}\mathbf{\hat{E}}_{j}})}\\
\hat{A}_{j}^{i} & =\frac{h_{i}A_{j}^{i}}{h_{j}} &  & \text{(no sum on \ensuremath{i} or \ensuremath{j}, with basis \ensuremath{\mathbf{\hat{E}}_{i}\mathbf{\hat{E}}^{j}})}
\end{aligned}
\end{equation}

$\bullet$ On generalizing the above pattern, the physical components
of a tensor of type ($m,n$) in a general orthogonal coordinate system
are given by:
\begin{equation}
\hat{A}_{b_{1}\ldots b_{n}}^{a_{1}\ldots a_{m}}=\frac{h_{a_{1}}\ldots h_{a_{m}}}{h_{b_{1}}\ldots h_{b_{n}}}A_{b_{1}\ldots b_{n}}^{a_{1}\ldots a_{m}}
\end{equation}

$\bullet$ As a consequence of the last points, in a space with a
well defined metric any tensor can be expressed in covariant or contravariant
(including mixed) or physical forms using different sets of basis
vectors. Moreover, these forms can be transformed from each other
using the raising and lowering operators and scale factors. As before,
for the Cartesian rectangular systems the covariant, contravariant
and physical components are the same where the Kronecker delta is
the metric tensor.

$\bullet$ For orthogonal coordinate systems, the two sets of normalized
covariant and contravariant basis vectors are identical as established
earlier, and hence the physical components related to the covariant
and contravariant components are identical as well. Consequently,
for orthogonal systems with orthonormal basis vectors, the covariant,
contravariant and physical components are identical.

$\bullet$ The physical components of a tensor may be represented
by the symbol of the tensor with subscripts denoting the coordinates
of the employed coordinate system. For instance, if $\mathbf{A}$
is a vector in a 3D space with contravariant components $A^{i}$ or
covariant components $A_{i}$, its physical components in Cartesian,
cylindrical, spherical and general curvilinear systems may be denoted
by ($A_{x},A_{y},A_{z}$), ($A_{\rho},A_{\phi},A_{z}$), ($A_{r},A_{\theta},A_{\phi}$)
and ($A_{u},A_{v},A_{w}$) respectively.

$\bullet$ For consistency and dimensional homogeneity, the tensors
in scientific applications are normally represented by their physical
components with a set of normalized unit base vectors. The invariance
of the tensor form then guarantees that the same tensor formulation
is valid regardless of any particular coordinate system where standard
tensor transformations can be used to convert from one form to another
without affecting the validity and invariance of the formulation.

\pagebreak{}

\section{Special Tensors}

$\bullet$ The subject of investigation of this section is those tensors
that form an essential part of the tensor calculus theory, namely
the Kronecker, the permutation and the metric tensors.

\subsection{Kronecker Tensor\label{subKronecker}}

$\bullet$ This is a rank-2 symmetric, constant, isotropic tensor
in all dimensions.

$\bullet$ It is defined as:
\begin{equation}
\delta_{ij}=\delta^{ij}=\delta_{\,\,j}^{i}=\delta_{i}^{\,\,j}\begin{cases}
1 & (i=j)\\
0\,\,\,\,\,\,\,\,\,\,\,\,\,\, & (i\neq j)
\end{cases}\label{eqKroneckerDefinitionNormal}
\end{equation}

$\bullet$ The generalized Kronecker delta is defined as:
\begin{equation}
\delta_{j_{1}\ldots j_{n}}^{i_{1}\ldots i_{n}}=\begin{cases}
\,\,\,\,\,1 & \left[(j_{1}\ldots j_{n})\text{ is even permutation of (\ensuremath{i_{1}\ldots i_{n})}}\right]\\
-1 & \left[(j_{1}\ldots j_{n})\text{ is odd permutation of (\ensuremath{i_{1}\ldots i_{n})}}\right]\\
\,\,\,\,\,0\,\,\,\,\,\,\,\,\,\,\,\,\,\, & \left[\text{repeated\,\,\ensuremath{j}'s}\right]
\end{cases}\label{eqGeneralizedKronecker}
\end{equation}

It can also be defined by the following $n\times n$ determinant:
\begin{equation}
\delta_{j_{1}\ldots j_{n}}^{i_{1}\ldots i_{n}}=\begin{vmatrix}\begin{array}{cccc}
\delta_{j_{1}}^{i_{1}} & \delta_{j_{2}}^{i_{1}} & \cdots & \delta_{j_{n}}^{i_{1}}\\
\delta_{j_{1}}^{i_{2}} & \delta_{j_{2}}^{i_{2}} & \cdots & \delta_{j_{n}}^{i_{2}}\\
\vdots & \vdots & \ddots & \vdots\\
\delta_{j_{1}}^{i_{n}} & \delta_{j_{2}}^{i_{n}} & \cdots & \delta_{j_{n}}^{i_{n}}
\end{array}\end{vmatrix}\label{eqGeneralizedKronecker2}
\end{equation}
where the $\delta_{j}^{i}$ entries in the determinant are the normal
Kronecker deltas as defined by Eq. \ref{eqKroneckerDefinitionNormal}.

$\bullet$ The relation between the rank-$n$ permutation tensor and
the generalized Kronecker delta in an $n$D space is given by:
\begin{equation}
\epsilon_{i_{1}i_{2}\ldots i_{n}}=\delta_{i_{1}i_{2}\ldots i_{n}}^{1\,2\ldots n}\,\,\,\,\,\,\,\,\,\,\,\,\,\,\,\,\,\,\,\,\,\,\,\&\,\,\,\,\,\,\,\,\,\,\,\,\,\,\,\,\,\,\,\,\,\,\,\epsilon^{i_{1}i_{2}\ldots i_{n}}=\delta_{1\,2\ldots n}^{i_{1}i_{2}\ldots i_{n}}
\end{equation}
Hence, the permutation tensor $\epsilon$ may be considered as a special
case of the generalized Kronecker delta. Consequently the permutation
tensor can be written as an $n\times n$ determinant consisting of
the normal Kronecker deltas.

$\bullet$ If we define
\begin{equation}
\delta_{lm}^{ij}=\delta_{lmk}^{ijk}
\end{equation}
then the well known $\epsilon-\delta$ relation (Eq. \ref{EqEpsilonDelta})
will take the following form:
\begin{equation}
\delta_{lm}^{ij}=\delta_{l}^{i}\delta_{m}^{j}-\delta_{m}^{i}\delta_{l}^{j}\label{eqDelDel}
\end{equation}
Other identities involving $\delta$ and $\epsilon$ can also be formulated
in terms of the generalized Kronecker delta.

\subsection{Permutation Tensor\label{subPermutation}}

$\bullet$ This tensor has a rank equal to the number of dimensions
of the space. Hence, a rank-$n$ permutation tensor has $n^{n}$ components.

$\bullet$ It is a relative tensor of weight $-1$ for its covariant
form and $+1$ for its contravariant form.

$\bullet$ It is isotropic and totally anti-symmetric in each pair
of its indices, i.e. it changes sign on swapping any two of its indices.

$\bullet$ It is a pseudo tensor since it acquires a minus sign under
improper orthogonal transformation of coordinates.

$\bullet$ The rank-$n$ permutation tensor is defined as:
\begin{equation}
\epsilon^{i_{1}i_{2}\ldots i_{n}}=\epsilon_{i_{1}i_{2}\ldots i_{n}}=\begin{cases}
\,\,\,\,\,1 & \left[\left(i_{1},i_{2},\ldots,i_{n}\right)\text{ is even permutation of (\ensuremath{1,2,\ldots,n})}\right]\\
-1 & \left[\left(i_{1},i_{2},\ldots,i_{n}\right)\text{ is odd permutation of (\ensuremath{1,2,\ldots,n})}\right]\\
\,\,\,\,\,0\,\,\,\,\,\,\,\,\,\,\,\,\,\, & \text{\ensuremath{\left[\mathrm{repeated\,index}\right]}}
\end{cases}\label{eqEpsilonnDefinition}
\end{equation}

$\bullet$ For the rank-$n$ permutation tensor we have:\footnote{The following formulae also apply to the contravariant form.}
\begin{equation}
\epsilon_{a_{1}a_{2}\cdots a_{n}}=\prod_{i=1}^{n-1}\left[\frac{1}{i!}\prod_{j=i+1}^{n}\left(a_{j}-a_{i}\right)\right]=\frac{1}{S(n-1)}\prod_{1\le i<j\le n}\left(a_{j}-a_{i}\right)\label{eqEpsilonSigma1}
\end{equation}
where $S(n-1)$ is the super-factorial function of $(n-1)$ which
is defined by:
\begin{equation}
S(k)=\prod_{i=1}^{k}i!=1!\cdot2!\cdot\ldots\cdot k!
\end{equation}

$\bullet$ A simpler formula for the rank-$n$ permutation tensor
can be obtained from the previous one by ignoring the magnitude of
the multiplication factors and taking their signs only, that is:
\begin{equation}
\epsilon_{a_{1}a_{2}\cdots a_{n}}=\prod_{1\le i<j\le n}\sigma\left(a_{j}-a_{i}\right)\label{eqEpsilonSigma}
\end{equation}
where
\begin{equation}
\sigma(k)=\begin{cases}
+1 & (k>0)\\
-1 & (k<0)\\
\,\,\,\,\,0\,\,\,\,\,\,\,\,\,\,\,\,\,\, & (k=0)
\end{cases}
\end{equation}

$\bullet$ The sign function in the previous point can be expressed
in a more direct form by dividing each argument of the multiplicative
factors in Eq. \ref{eqEpsilonSigma} by its absolute value, noting
that none of these factors is zero, and hence Eq. \ref{eqEpsilonSigma}
becomes:
\begin{equation}
\epsilon_{a_{1}a_{2}\cdots a_{n}}=\prod_{1\le i<j\le n}\frac{\left(a_{j}-a_{i}\right)}{\left|a_{j}-a_{i}\right|}
\end{equation}

$\bullet$ For the rank-3 permutation tensor we have:
\begin{equation}
\epsilon^{ijk}\epsilon_{lmn}=\begin{vmatrix}\begin{array}{ccc}
\delta_{l}^{i} & \delta_{m}^{i} & \delta_{n}^{i}\\
\delta_{l}^{j} & \delta_{m}^{j} & \delta_{n}^{j}\\
\delta_{l}^{k} & \delta_{m}^{k} & \delta_{n}^{k}
\end{array}\end{vmatrix}
\end{equation}
\begin{equation}
\epsilon^{ijk}\epsilon_{lmk}=\delta_{l}^{i}\delta_{m}^{j}-\delta_{m}^{i}\delta_{l}^{j}\label{EqEpsilonDelta}
\end{equation}

$\bullet$ For the rank-$n$ permutation tensor we have:
\begin{equation}
\epsilon^{i_{1}i_{2}\cdots i_{n}}\,\epsilon_{i_{1}i_{2}\cdots i_{n}}=n!
\end{equation}

$\bullet$ For the rank-$n$ permutation tensor we have:
\begin{equation}
\epsilon^{i_{1}i_{2}\cdots i_{n}}\,\epsilon_{j_{1}j_{2}\cdots j_{n}}=\begin{vmatrix}\begin{array}{cccc}
\delta_{j_{1}}^{i_{1}} & \delta_{j_{2}}^{i_{1}} & \cdots & \delta_{j_{n}}^{i_{1}}\\
\delta_{j_{1}}^{i_{2}} & \delta_{j_{2}}^{i_{2}} & \cdots & \delta_{j_{n}}^{i_{2}}\\
\vdots & \vdots & \ddots & \vdots\\
\delta_{j_{1}}^{i_{n}} & \delta_{j_{2}}^{i_{n}} & \cdots & \delta_{j_{n}}^{i_{n}}
\end{array}\end{vmatrix}\label{eqEpsilon2}
\end{equation}

$\bullet$ On comparing Eqs. \ref{eqGeneralizedKronecker2} and \ref{eqEpsilon2}
we obtain the following identity:
\begin{equation}
\delta_{j_{1}\ldots j_{n}}^{i_{1}\ldots i_{n}}=\epsilon^{i_{1}\ldots i_{n}}\,\epsilon_{j_{1}\ldots j_{n}}
\end{equation}

$\bullet$ Based on the previous point, the generalized Kronecker
delta is the result of multiplying two relative tensors one of weight
$w=+1$ and the other of weight $w=-1$ and hence the generalized
Kronecker delta has a weight of $w=0$; therefore the generalized
Kronecker delta is an absolute tensor.\footnote{The multiplication of relative tensors produces a tensor whose weight
is the sum of the weights of the original tensors.}

$\bullet$ As has been stated previously, $\epsilon^{i_{1}\ldots i_{n}}$
and $\epsilon_{i_{1}\ldots i_{n}}$ are relative tensors of weight
$+1$ and $-1$ respectively. It is desirable to define absolute covariant
and contravariant forms of the permutation tensor, marked with underline,
by the following relations:\footnote{The contravariant form requires a sign function with details out of
scope of the present text (see \cite{ZwillingerBook2012}); however,
for the rank-3 permutation tensor which is the one used mostly in
the forthcoming notes the above expression stands as it is.}
\begin{equation}
\underline{\epsilon}_{i_{1}\ldots i_{n}}=\sqrt{g}\,\epsilon_{i_{1}\ldots i_{n}}\,\,\,\,\,\,\,\,\,\,\,\,\,\,\,\,\,\,\,\&\,\,\,\,\,\,\,\,\,\,\,\,\,\,\,\,\,\underline{\epsilon}^{i_{1}\ldots i_{n}}=\frac{1}{\sqrt{g}}\epsilon^{i_{1}\ldots i_{n}}\label{eqepsilonunderline}
\end{equation}
where $g$ is the determinant of the covariant metric tensor $g_{pq}$.

$\bullet$ The $\epsilon-\delta$ identity (Eqs. \ref{eqDelDel} and
\ref{EqEpsilonDelta}) can be generalized by employing the metric
tensor with the absolute permutation tensor:
\begin{equation}
g^{ij}\underline{\epsilon}_{ikl}\underline{\epsilon}_{jmn}=g_{km}g_{ln}-g_{kn}g_{lm}
\end{equation}

\subsection{Metric Tensor\label{subMetricTensor}}

$\bullet$ One of the main objectives of the metric, which is a rank-2
symmetric absolute non-singular\footnote{We mean that the matrix representing the tensor is invertible and
hence its determinant does not vanish at any point of the space.} tensor, is to generalize the concept of distance to general curvilinear
coordinate frames and hence maintain the invariance of distance in
different coordinate systems. This tensor is also used to raise and
lower indices and thus facilitate the transformation between the covariant
and contravariant types.

$\bullet$ In general, the coordinate system and the space metric
are independent entities. Yes, some coordinate systems may be defined
by having a specific metric in which case the two are correlated.
This is the case in the Cartesian coordinate systems which are based
in their definition on presuming an underlying Euclidean metric.

$\bullet$ The components of the metric tensor are given by:
\begin{equation}
g_{ij}=\mathbf{E}_{i}\cdot\mathbf{E}_{j}\,\,\,\,\,\,\,\,\,\,\,\,\,\,\,\,\,\&\,\,\,\,\,\,\,\,\,\,\,\,\,\,\,\,\,\,g^{ij}=\mathbf{E}^{i}\cdot\mathbf{E}^{j}\label{eqEg}
\end{equation}
where the indexed $\mathbf{E}$ are the covariant and contravariant
basis vectors as defined previously in $\S$ \ref{subBasisVectors}.
Because of these relations, the vectors $\mathbf{E}_{i}$ and $\mathbf{E}^{i}$
may be denoted by $\mathbf{g}_{i}$ and $\mathbf{g}^{i}$ respectively
which is more suggestive of their relation to the metric tensor.

$\bullet$ As a consequence of the last point, the covariant metric
tensor can also be defined as:
\begin{equation}
g_{ij}=\frac{\partial x^{k}}{\partial u^{i}}\frac{\partial x^{k}}{\partial u^{j}}\label{eqgxu}
\end{equation}
where
\begin{equation}
x^{k}=x^{k}\left(u^{1},\ldots,u^{n}\right)\,\,\,\,\,\,\,\,\,\,\,\,\,\,\,\,\,\,\text{(\ensuremath{k=1,\ldots,n})}
\end{equation}
are independent coordinates in an $n$D space with a rectangular Cartesian
system, and $u^{i}$ $\text{(\ensuremath{i=1,\ldots,n})}$ are independent
generalized curvilinear coordinates. Similarly for the contravariant
metric tensor we have:
\begin{equation}
g^{ij}=\frac{\partial u^{i}}{\partial x^{k}}\frac{\partial u^{j}}{\partial x^{k}}\label{eqgux}
\end{equation}

$\bullet$ The coefficients of the metric tensor may also be considered
as the components of the unit tensor in its two variance forms, that
is:
\begin{equation}
\boldsymbol{\delta}=g_{ij}\mathbf{E}^{i}\mathbf{E}^{j}=g^{ij}\mathbf{E}_{i}\mathbf{E}_{j}
\end{equation}

$\bullet$ As stated already, the basis vectors, whether covariant
or contravariant, in general coordinate systems are not necessarily
mutually orthogonal and hence the metric tensor is not diagonal in
general since the dot products given in Eqs. \ref{eqEg}, \ref{eqgxu}
and \ref{eqgux} are not necessarily zero when $i\ne j$. Moreover,
since those basis vectors are not necessarily of unit length, the
entries of the metric tensor are not of unit magnitude in general.
However, since the dot product of vectors is a commutative operation,
the metric tensor is necessarily symmetric.

$\bullet$ The entries of the metric tensor, including the diagonal
elements, can be positive or negative.

$\bullet$ The covariant and contravariant forms of the metric tensor
are inverses of each other and hence:
\begin{equation}
g^{ik}g_{kj}=\delta_{\,\,j}^{i}\,\,\,\,\,\,\,\,\,\,\,\,\,\,\,\,\&\,\,\,\,\,\,\,\,\,\,\,\,\,\,\,\,g_{ik}g^{kj}=\delta_{i}^{\,\,j}
\end{equation}
where these equations can be seen as a matrix multiplication (row$\times$column).

$\bullet$ A result from the previous points is that:
\begin{equation}
\begin{aligned}\left(\mathbf{E}^{i}\cdot\mathbf{E}^{j}\right)\left(\mathbf{E}_{j}\cdot\mathbf{E}_{k}\right) & =g^{ij}g_{jk}=\delta_{\,\,k}^{i}\\
\left(\mathbf{E}_{i}\cdot\mathbf{E}_{j}\right)\left(\mathbf{E}^{j}\cdot\mathbf{E}^{k}\right) & =g_{ij}g^{jk}=\delta_{i}^{\,\,k}
\end{aligned}
\end{equation}

$\bullet$ As the metric tensor has an inverse, it should not be singular
and hence its determinant, which is in general a function of coordinates
like the metric tensor itself, should not vanish at any point in the
space, that is:
\begin{equation}
g(u^{1},\ldots,u^{n})=\mathrm{det}\left(g_{ij}\right)\ne0
\end{equation}

$\bullet$ The mixed type metric tensor is given by:
\begin{equation}
g_{\,\,j}^{i}=\mathbf{E}^{i}\cdot\mathbf{E}_{j}=\delta_{\,\,j}^{i}\,\,\,\,\,\,\,\,\,\,\,\,\,\,\,\&\,\,\,\,\,\,\,\,\,\,\,\,\,\,\,\,g_{i}^{\,\,j}=\mathbf{E}_{i}\cdot\mathbf{E}^{j}=\delta_{i}^{\,\,j}\label{eqEiEjMix}
\end{equation}
and hence it is the identity tensor. These equations represent the
fact that the covariant and contravariant basis vectors are reciprocal
sets.

$\bullet$ From the previous points, it can be concluded that the
metric tensor is in fact a transformation of the Kronecker delta in
its different variance types from a rectangular system to a general
curvilinear system, that is:
\begin{equation}
\begin{aligned}g_{ij} & =\frac{\partial x^{k}}{\partial u^{i}}\frac{\partial x^{l}}{\partial u^{j}}\delta_{kl}=\frac{\partial x^{k}}{\partial u^{i}}\frac{\partial x^{k}}{\partial u^{j}}=\mathbf{E}_{i}\cdot\mathbf{E}_{j} & \,\,\,\,\,\,\,\,\, & \text{(covariant)}\\
g^{ij} & =\frac{\partial u^{i}}{\partial x^{k}}\frac{\partial u^{j}}{\partial x^{l}}\delta^{kl}=\frac{\partial u^{i}}{\partial x^{k}}\frac{\partial u^{j}}{\partial x^{k}}=\mathbf{E}^{i}\cdot\mathbf{E}^{j} &  & \text{(contravariant)}\\
g_{j}^{i} & =\frac{\partial u^{i}}{\partial x^{k}}\frac{\partial x^{l}}{\partial u^{j}}\delta_{l}^{k}=\frac{\partial u^{i}}{\partial x^{k}}\frac{\partial x^{k}}{\partial u^{j}}=\mathbf{E}^{i}\cdot\mathbf{E}_{j} &  & \text{(mixed)}
\end{aligned}
\end{equation}

$\bullet$ Because of the relations:
\begin{equation}
\begin{aligned}A^{i} & =\mathbf{A}\cdot\mathbf{E}^{i}=A_{j}\mathbf{E}^{j}\cdot\mathbf{E}^{i}=A_{j}g^{ji}\\
A_{i} & =\mathbf{A}\cdot\mathbf{E}_{i}=A^{j}\mathbf{E}_{j}\cdot\mathbf{E}_{i}=A^{j}g_{ji}
\end{aligned}
\end{equation}
the metric tensor is used as an operator for raising and lowering
indices and hence facilitating the transformation between the covariant
and contravariant types of vectors. By a similar argument, the above
can be easily generalized where the contravariant metric tensor is
used for raising covariant indices and the covariant metric tensor
is used for lowering contravariant indices of tensors of any rank,
e.g.
\begin{equation}
A_{\,\,k}^{i}=g^{ij}A_{jk}\,\,\,\,\,\,\,\,\,\,\,\,\,\,\,\,\,\,\,\&\,\,\,\,\,\,\,\,\,\,\,\,\,\,\,\,\,A_{i}^{\,\,kl}=g_{ij}A^{jkl}
\end{equation}
Consequently, any tensor in a Riemannian space with well-defined metric
can be cast into covariant or contravariant or mixed forms.\footnote{For mixed form the rank should be $>1$.}

$\bullet$ In the raising and lowering of index operations the metric
tensor acts, like a Kronecker delta, as an index replacement operator
as well as shifting the index position.

$\bullet$ In general, the order of the raised and lowered indices
is important and hence
\begin{equation}
g^{ik}A_{jk}=A_{j}^{\,\,\,i}\,\,\,\,\,\,\,\,\,\,\,\,\,\,\mathrm{and}\,\,\,\,\,\,\,\,\,\,\,\,\,\,\,\,g^{ik}A_{kj}=A_{\,\,\,j}^{i}
\end{equation}
are different unless the tensor is symmetric in its two indices, i.e.
$A_{jk}=A_{kj}$. A dot may be used to indicate the original position
of the shifted index and hence the order of the indices is recorded,
e.g. $A_{j\,\cdot}^{\,\,\,i}$ and $A_{\cdot\,j}^{i}$ for the above
examples respectively, although this is redundant in the case of symmetry.\footnote{Dots may also be inserted in the tensor symbols to remove any ambiguity
about the order of the indices even without the action of the raising
and lowering operators.}

$\bullet$ Raising and lowering of indices is a reversible process;
hence keeping a record of the original position of the shifted indices
will facilitate the reversal.

$\bullet$ For a space with a coordinate system in which the metric
tensor can be cast into a diagonal form with all the diagonal entries
being of unity magnitude (i.e. $\pm1$) the metric is called flat.

$\bullet$ If $g$ and $\bar{g}$ are the determinants of the covariant
metric tensor in the unbarred and barred systems respectively, i.e.
$g=\mathrm{det}\left(g_{ij}\right)$ and $\bar{g}=\mathrm{det}\left(\bar{g}_{ij}\right)$,
then
\begin{equation}
\bar{g}=J^{2}g\,\,\,\,\,\,\,\,\,\,\,\,\,\,\,\,\,\&\,\,\,\,\,\,\,\,\,\,\,\,\,\,\,\,\,\sqrt{\bar{g}}=J\sqrt{g}\label{eqgJ}
\end{equation}
where $J$ ($=\left|\frac{\partial u}{\partial\bar{u}}\right|$) is
the Jacobian of the transformation between the unbarred and barred
systems. Consequently, the determinant of the covariant metric and
its square root are relative scalar invariants of weight $+2$ and
$+1$ respectively.

$\bullet$ A ``conjugate'' or ``associated'' tensor of a tensor
in a metric space is a tensor obtained by inner product multiplication,
once or more, of the original tensor by the covariant or contravariant
forms of the metric tensor.

$\bullet$ All tensors associated with a particular tensor through
the metric tensor represent the same tensor but in different reference
frames since the association is no more than raising or lowering indices
by the metric tensor which is equivalent to a representation of the
components of the tensor relative to different basis sets.

$\bullet$ A sufficient and necessary condition for the components
of the metric tensor to be constants in a given coordinate system
is that the Christoffel symbols of the first or second kind vanish
identically (refer to \ref{subChristoffelSymbols}).

$\bullet$ The metric tensor behaves as a constant with respect to
covariant and absolute differentiation (see $\S$ \ref{subCovariantDerivative}
and $\S$ \ref{subAbsoluteDerivative}). Hence, in all coordinate
systems the covariant and absolute derivatives of the metric tensor
are zero; moreover, the covariant and absolute derivative operators
bypass the metric tensor in differentiating inner and outer products
of tensors involving the metric tensor.

$\bullet$ In general orthogonal coordinate systems in $n$D spaces
the metric tensor and its inverse are diagonal, that is:
\begin{equation}
g_{ij}=g^{ij}=0\,\,\,\,\,\,\,\,\,\,\,\,\,\,\,\,\,\text{(\ensuremath{i\ne j})}
\end{equation}
moreover, we have:
\begin{equation}
g_{ii}=\left(h_{i}\right)^{2}=\frac{1}{g^{ii}}\,\,\,\,\,\,\,\,\,\,\,\,\,\,\,\,\,\,\,\,\,\,\text{(no sum on \ensuremath{i})}
\end{equation}
\begin{equation}
\mathrm{det}\left(g_{ij}\right)=g=g_{11}g_{22}\ldots g_{nn}=\prod_{i}\left(h_{i}\right)^{2}
\end{equation}
\begin{equation}
\mathrm{det}\left(g^{ij}\right)=\frac{1}{g}=\frac{1}{g_{11}g_{22}\ldots g_{nn}}=\left[\prod_{i}\left(h_{i}\right)^{2}\right]^{-1}
\end{equation}
where $h_{i}$ ($=\left|\mathbf{E}_{i}\right|$) are the scale factors,
as described previously.

$\bullet$ A Riemannian metric, $g_{ij}$, in a particular coordinate
system is a Euclidean metric if it can be transformed to the identity
tensor, $\delta_{ij}$, by a permissible coordinate transformation.

$\bullet$ The Minkowski metric, which is the metric tensor of special
relativity, is given by one of the following two forms:
\begin{equation}
\left[g_{ij}\right]=\left[g^{ij}\right]=\left[\begin{array}{cccc}
1 & 0 & 0 & 0\\
0 & -1 & 0 & 0\\
0 & 0 & -1 & 0\\
0 & 0 & 0 & -1
\end{array}\right]\,\,\,\,\,\,\,\,\,\,\,\,\,\,\,\,\,\,\,\,\,\,\,\,\,\,\,\,\left[g_{ij}\right]=\left[g^{ij}\right]=\left[\begin{array}{cccc}
-1 & 0 & 0 & 0\\
0 & 1 & 0 & 0\\
0 & 0 & 1 & 0\\
0 & 0 & 0 & 1
\end{array}\right]
\end{equation}
Consequently, the line element $ds$ can be imaginary.

$\bullet$ The partial derivatives of the covariant and contravariant
metric tensors satisfy the following identities:
\begin{equation}
\begin{aligned}\partial_{k}g_{ij} & =-g_{mj}g_{ni}\partial_{k}g^{nm}\\
\partial_{k}g^{ij} & =-g^{mj}g^{in}\partial_{k}g_{nm}
\end{aligned}
\end{equation}

$\bullet$ In the following subsections, we investigate a number of
mathematical objects whose definitions and applications are dependent
on the metric tensor.

\subsubsection{Dot Product}

$\bullet$ The dot product of two basis vectors in general curvilinear
coordinates was given earlier in this section. This will be used in
the following points to develop expressions for the dot product of
vectors and tensors in general.

$\bullet$ The dot product of two vectors, $\mathbf{A}$ and $\mathbf{B}$,
in general curvilinear coordinates using their covariant and contravariant
forms, as well as opposite forms, is given by:
\begin{eqnarray}
\mathbf{A}\cdot\mathbf{B} & = & A_{i}\mathbf{E}^{i}\cdot B_{j}\mathbf{E}^{j}=A_{i}B_{j}\mathbf{E}^{i}\cdot\mathbf{E}^{j}=g^{ij}A_{i}B_{j}=A^{j}B_{j}=A_{i}B^{i}\nonumber \\
\mathbf{A}\cdot\mathbf{B} & = & A^{i}\mathbf{E}_{i}\cdot B^{j}\mathbf{E}_{j}=A^{i}B^{j}\mathbf{E}_{i}\cdot\mathbf{E}_{j}=g_{ij}A^{i}B^{j}=A_{j}B^{j}=A^{i}B_{i}\\
\mathbf{A}\cdot\mathbf{B} & = & A_{i}\mathbf{E}^{i}\cdot B^{j}\mathbf{E}_{j}=A_{i}B^{j}\mathbf{E}^{i}\cdot\mathbf{E}_{j}=\delta_{\,j}^{i}A_{i}B^{j}=A_{j}B^{j}\nonumber \\
\mathbf{A}\cdot\mathbf{B} & = & A^{i}\mathbf{E}_{i}\cdot B_{j}\mathbf{E}^{j}=A^{i}B_{j}\mathbf{E}_{i}\cdot\mathbf{E}^{j}=\delta_{i}^{\,j}A^{i}B_{j}=A^{i}B_{i}\nonumber
\end{eqnarray}
In brief, the dot product of two vectors is the dot product of their
two basis vectors multiplied algebraically by the algebraic product
of their components. Because the dot product of basis vectors is a
metric tensor, the metric tensor will act on the components by raising
or lowering the index of one component or by replacing the index of
a component.

$\bullet$ The dot product operations outlined in the previous point
can be easily extended to tensors of higher ranks where the covariant
and contravariant forms of the components and basis vectors are treated
in a similar manner to the above to obtain the dot product. For instance,
the dot product of a rank-2 tensor of contravariant components $A^{ij}$
and a vector of covariant components $B_{k}$ is given by:
\begin{equation}
\mathbf{A}\cdot\mathbf{B}=\left(A^{ij}\mathbf{E}_{i}\mathbf{E}_{j}\right)\cdot\left(B_{k}\mathbf{E}^{k}\right)=A^{ij}B_{k}\left(\mathbf{E}_{i}\mathbf{E}_{j}\cdot\mathbf{E}^{k}\right)=A^{ij}B_{k}\mathbf{E}_{i}\delta_{j}^{\,k}=A^{ij}B_{j}\mathbf{E}_{i}
\end{equation}
that is, the $i^{th}$ component of this product, which is a contravariant
vector, is:
\begin{equation}
\left[\mathbf{A}\cdot\mathbf{B}\right]^{i}=A^{ij}B_{j}
\end{equation}

$\bullet$ From the previous points, the dot product in general curvilinear
coordinates occurs between two vectors of opposite variance type.
Therefore, to obtain the dot product of two vectors of the same variance
type, one of the vectors should be converted to the opposite type
by the raising/lowering operator, followed by the inner product operation.
This can be generalized to the dot product of higher-rank tensors
where the two contracted indices of the dot product should be of opposite
variance type and hence the index-shifting operator in the form of
the metric tensor should be used, if necessary, to achieve this.

$\bullet$ The generalized dot product of two tensors is an invariant
under permissible coordinate transformations.

\subsubsection{Cross Product}

$\bullet$ The cross product of two covariant basis vectors in general
curvilinear coordinates is given by:
\begin{equation}
\mathbf{E}_{i}\times\mathbf{E}_{j}=\frac{\partial x^{l}}{\partial u^{i}}\mathbf{e}_{l}\times\frac{\partial x^{m}}{\partial u^{j}}\mathbf{e}_{m}=\frac{\partial x^{l}}{\partial u^{i}}\frac{\partial x^{m}}{\partial u^{j}}\mathbf{e}_{l}\times\mathbf{e}_{m}=\frac{\partial x^{l}}{\partial u^{i}}\frac{\partial x^{m}}{\partial u^{j}}\epsilon_{lmn}\mathbf{e}_{n}
\end{equation}
where the indexed $x$ and $u$ are the coordinates of Cartesian and
general curvilinear systems respectively, the indexed $\mathbf{e}$
are the Cartesian base vectors\footnote{For Cartesian systems, there is no difference between covariant and
contravariant tensors and hence $\mathbf{e}_{i}=\mathbf{e}^{i}$.
We also note that for Cartesian systems $g=1$.} and $\epsilon_{lmn}=\epsilon^{lmn}$ is the permutation relative
tensor as defined in Eq. \ref{eqEpsilonnDefinition}. Now since $\mathbf{e}_{n}=\mathbf{e}^{n}=\frac{\partial x^{n}}{\partial u^{k}}\mathbf{E}^{k}$,
the last equation becomes:
\begin{equation}
\mathbf{E}_{i}\times\mathbf{E}_{j}=\frac{\partial x^{l}}{\partial u^{i}}\frac{\partial x^{m}}{\partial u^{j}}\frac{\partial x^{n}}{\partial u^{k}}\epsilon_{lmn}\mathbf{E}^{k}=\underline{\epsilon}_{ijk}\mathbf{E}^{k}
\end{equation}
where the underlined absolute covariant permutation tensor is defined
as:
\begin{equation}
\underline{\epsilon}_{ijk}=\frac{\partial x^{l}}{\partial u^{i}}\frac{\partial x^{m}}{\partial u^{j}}\frac{\partial x^{n}}{\partial u^{k}}\epsilon_{lmn}\label{eqepsilonAbsoluteCo}
\end{equation}
So the final result is:
\begin{equation}
\mathbf{E}_{i}\times\mathbf{E}_{j}=\underline{\epsilon}_{ijk}\mathbf{E}^{k}\label{eqEiEjCov}
\end{equation}
By a similar reasoning, we obtain the following expression for the
cross product of two contravariant basis vectors in general curvilinear
coordinates:
\begin{equation}
\mathbf{E}^{i}\times\mathbf{E}^{j}=\underline{\epsilon}^{ijk}\mathbf{E}_{k}\label{eqEiEjCon}
\end{equation}
where the absolute contravariant permutation tensor is defined by:
\begin{equation}
\underline{\epsilon}^{ijk}=\frac{\partial u^{i}}{\partial x^{l}}\frac{\partial u^{j}}{\partial x^{m}}\frac{\partial u^{k}}{\partial x^{n}}\epsilon^{lmn}\label{eqepsilonunderlineJacobian}
\end{equation}

$\bullet$ Considering Eq. \ref{eqepsilonunderline}, the above equations
can also be expressed as:
\begin{equation}
\mathbf{E}_{i}\times\mathbf{E}_{j}=\underline{\epsilon}_{ijk}\mathbf{E}^{k}=\sqrt{g}\epsilon_{ijk}\mathbf{E}^{k}\label{eqECoCrossPro}
\end{equation}
\begin{equation}
\mathbf{E}^{i}\times\mathbf{E}^{j}=\underline{\epsilon}^{ijk}\mathbf{E}_{k}=\frac{\epsilon^{ijk}}{\sqrt{g}}\mathbf{E}_{k}\label{eqECovaCrossPro}
\end{equation}
where $\epsilon_{ijk}=\epsilon^{ijk}$ are as defined previously (Eq.
\ref{eqEpsilonnDefinition}).

$\bullet$ The cross product of non-basis vectors follows similar
rules to those outlined above for the basis vectors; the only difference
is that the algebraic product of the components is used as a scale
factor for the cross product of their basis vectors. For example,
the cross product of two contravariant vectors, $A^{i}$ and $B^{j}$,
is given by:
\begin{equation}
\mathbf{A}\times\mathbf{B}=\left(A^{i}\mathbf{E}_{i}\right)\times\left(B^{j}\mathbf{E}_{j}\right)=A^{i}B^{j}\left(\mathbf{E}_{i}\times\mathbf{E}_{j}\right)=\underline{\epsilon}_{ijk}A^{i}B^{j}\mathbf{E}^{k}
\end{equation}
that is, the $k^{th}$ component of this product, which is a vector
with covariant components, is:
\begin{equation}
\left[\mathbf{A}\times\mathbf{B}\right]_{k}=\underline{\epsilon}_{ijk}A^{i}B^{j}
\end{equation}
Similarly, the cross product of two covariant vectors, $A_{i}$ and
$B_{j}$, is given by:
\begin{equation}
\mathbf{A}\times\mathbf{B}=\left(A_{i}\mathbf{E}^{i}\right)\times\left(B_{j}\mathbf{E}^{j}\right)=A_{i}B_{j}\left(\mathbf{E}^{i}\times\mathbf{E}^{j}\right)=\underline{\epsilon}^{ijk}A_{i}B_{j}\mathbf{E}_{k}
\end{equation}
with the $k^{th}$ contravariant component being given by:
\begin{equation}
\left[\mathbf{A}\times\mathbf{B}\right]^{k}=\underline{\epsilon}^{ijk}A_{i}B_{j}
\end{equation}

\subsubsection{Line Element}

$\bullet$ The displacement differential vector in general curvilinear
coordinate systems is given by:
\begin{equation}
d\mathbf{r}=\frac{\partial\mathbf{r}}{\partial u^{i}}du^{i}=\mathbf{E}_{i}du^{i}=\sum_{i}\left|\mathbf{E}_{i}\right|\frac{\mathbf{E}_{i}}{\left|\mathbf{E}_{i}\right|}du^{i}=\sum_{i}\left|\mathbf{E}_{i}\right|\hat{\mathbf{E}}_{i}du^{i}
\end{equation}
where $\mathbf{r}$ is the position vector as defined previously.

$\bullet$ The line element $ds$, which may also be called the differential
of arc length, in general curvilinear coordinate systems is given
by:
\begin{equation}
\left(ds\right)^{2}=d\mathbf{r}\cdot d\mathbf{r}=\mathbf{E}_{i}du^{i}\cdot\mathbf{E}_{j}du^{j}=\left(\mathbf{E}_{i}\cdot\mathbf{E}_{j}\right)du^{i}du^{j}=g_{ij}du^{i}du^{j}\label{eqdsSquared}
\end{equation}
where $g_{ij}$ is the covariant metric tensor.

$\bullet$ For orthogonal coordinate systems, the metric tensor is
given by:
\begin{equation}
g_{ij}=\begin{cases}
0 & (i\ne j)\\
\left(h_{i}\right)^{2}\,\,\,\,\,\,\,\,\,\,\,\,\,\, & (i=j)
\end{cases}
\end{equation}
where $h_{i}$ is the scale factor of the respective coordinate $u^{i}$.
Hence, the last part of Eq. \ref{eqdsSquared} becomes:
\begin{equation}
\left(ds\right)^{2}=\sum_{i}\left(h_{i}\right)^{2}du^{i}du^{i}\label{eqds2unbarred}
\end{equation}
with no cross terms (i.e. terms of products involving more than one
coordinate like $du^{i}du^{j}$ where $i\ne j$) which are generally
present in the case of non-orthogonal curvilinear systems.

$\bullet$ On conducting a transformation from one coordinate system
to another coordinate system, marked with barred coordinates, $\bar{u}$,
the line element will be expressed in the new system as:
\begin{equation}
\left(ds\right)^{2}=\bar{g}_{ij}d\bar{u}^{i}d\bar{u}^{j}\label{eqds2barred}
\end{equation}
Since the line element is an invariant quantity, the same symbol $(ds)^{2}$
is used in both Eqs. \ref{eqdsSquared} and \ref{eqds2barred}.

\subsubsection{Surface Element}

$\bullet$ In general curvilinear coordinates of a 3D space, an infinitesimal
element of area on the surface $u^{1}=c_{1}$, where $c_{1}$ is a
constant, is obtained by taking the magnitude of the cross product
of the displacement vectors in the directions of the other two coordinates
on that surface. Hence, the generalized differential of area element
on the surface $u^{1}=c_{1}$ is given by:
\begin{equation}
\begin{aligned}dA(u^{1}=c_{1}) & =\left|d\mathbf{r}_{2}\times d\mathbf{r}_{3}\right| & \,\,\,\,\,\,\,\,\,\,\,\,\,\,\,\,\,\,\,\,\\
 & =\left|\mathbf{E}_{2}\times\mathbf{E}_{3}\right|du^{2}du^{3}\\
 & =\left|\underline{\epsilon}_{231}\mathbf{E}^{1}\right|du^{2}du^{3} &  & \text{(Eq. \ref{eqEiEjCov})}\\
 & =\left|\underline{\epsilon}_{231}\right|\left|\mathbf{E}^{1}\right|du^{2}du^{3}\\
 & =\sqrt{g}\sqrt{\mathbf{E}^{1}\cdot\mathbf{E}^{1}}\,du^{2}du^{3} &  & \text{(Eqs. \ref{eqepsilonunderline} \& \ref{eqVectorMagnitudeCon})}\\
 & =\sqrt{g}\,\sqrt{g^{11}}\,du^{2}du^{3} &  & \text{(Eq. \ref{eqEg})}\\
 & =\sqrt{gg^{11}}\,du^{2}du^{3}
\end{aligned}
\end{equation}

$\bullet$ On generalizing the above argument, the differential area
element in a 3D space on the surface $u^{i}=c_{i}$ ($i=1,2,3$) where
$c_{i}$ is a constant is given by:
\begin{equation}
dA(u^{i}=c_{i})=\sqrt{gg^{ii}}du^{j}du^{k}\,\,\,\,\,\,\,\,\,\,\,\,\,\,\,\,\,\,\,\,\,\,\text{(\ensuremath{i\ne j\ne k}, no sum on \ensuremath{i})}\label{eqdAgeneral}
\end{equation}

$\bullet$ In general orthogonal  coordinates in a 3D space we have:
\begin{equation}
\sqrt{gg^{ii}}=\sqrt{(h_{i})^{2}(h_{j})^{2}(h_{k})^{2}\frac{1}{(h_{i})^{2}}}=h_{j}h_{k}\,\,\,\,\,\,\,\,\,\,\,\,\,\,\,\,\,\,\text{(\ensuremath{i\ne j\ne k}, no sum on any index)}
\end{equation}
and hence Eq. \ref{eqdAgeneral} becomes:
\begin{equation}
dA(u^{i}=c_{i})=h_{j}h_{k}du^{j}du^{k}\,\,\,\,\,\,\,\,\,\,\,\,\,\,\,\,\,\,\,\,\,\,\text{(\ensuremath{i\ne j\ne k}, no sum on any index)}
\end{equation}
The last formula represents the area of a surface differential with
sides $h_{j}du^{j}$ and $h_{k}du^{k}$ (no sum on $j,k$).

\subsubsection{Volume Element}

$\bullet$ In general curvilinear coordinates of a 3D space, an infinitesimal
element of volume, represented by a parallelepiped spanned by the
three displacement vectors $d\mathbf{r}_{i}\,\,\,(i=1,2,3)$, is obtained
by taking the magnitude of the scalar triple product of these vectors.
Hence, the generalized differential volume element is given by:
\begin{equation}
\begin{aligned}dV & =\left|d\mathbf{r}_{1}\cdot\left(d\mathbf{r}_{2}\times d\mathbf{r}_{3}\right)\right| & \,\,\,\,\,\,\,\,\,\,\,\,\,\,\,\,\,\,\,\,\,\,\,\,\,\,\,\,\,\,\\
 & =\left|\mathbf{E}_{1}\cdot\left(\mathbf{E}_{2}\times\mathbf{E}_{3}\right)\right|du^{1}du^{2}du^{3}\\
 & =\left|\mathbf{E}_{1}\cdot\underline{\epsilon}_{231}\mathbf{E}^{1}\right|du^{1}du^{2}du^{3} &  & \text{(Eq. \ref{eqEiEjCov})}\\
 & =\left|\mathbf{E}_{1}\cdot\mathbf{E}^{1}\right|\left|\underline{\epsilon}_{231}\right|du^{1}du^{2}du^{3}\\
 & =\left|\delta_{1}^{\,1}\right|\left|\underline{\epsilon}_{231}\right|du^{1}du^{2}du^{3} &  & \text{(Eq. \ref{eqEiEjMix})}\\
 & =\sqrt{g}\,du^{1}du^{2}du^{3} &  & \text{(Eq. \ref{eqepsilonunderline})}\\
 & =J\,du^{1}du^{2}du^{3} &  & \text{(Eq. \ref{eqgJsquared})}
\end{aligned}
\end{equation}
where $g$ is the determinant of the covariant metric tensor $g_{ij}$,
and $J$ is the Jacobian\footnote{Due to the freedom of choice in the order of the variables, which
is related to the choice of the system handedness hence affecting
the sign of the determinant Jacobian, the sign of the determinant
should be adjusted if necessary to have a proper sign for the volume
element.} of the transformation as defined previously. The last line in the
last equation is particularly relevant to the case of change of variables
in multivariate integrals where the Jacobian facilitates the transformation.

$\bullet$ The formulae in the last point for a 3D space can be extended
to the differential of generalized volume element\footnote{Generalized volume elements are used, for instance, to represent the
change of variables in multi-variable integrations.} in general curvilinear coordinates in an $n$D space as follow:
\begin{equation}
dV=\sqrt{g}du^{1}\ldots du^{n}=J\,du^{1}\ldots du^{n}
\end{equation}

$\bullet$ In general orthogonal  coordinate systems in a 3D space,
the above formulae become:
\begin{equation}
dV=h_{1}h_{2}h_{3}\,du^{1}du^{2}du^{3}
\end{equation}
where $h_{1},h_{2}$ and $h_{3}$ are the scale factors. The last
formula represents the volume of a parallelepiped with edges $h_{1}du^{1}$,
$h_{2}du^{2}$ and $h_{3}du^{3}$.

\subsubsection{Magnitude of Vector}

$\bullet$ The magnitude of a contravariant vector $\mathbf{A}$ is
given by:
\begin{equation}
\left|\mathbf{A}\right|=\sqrt{\mathbf{A}\cdot\mathbf{A}}=\sqrt{\left(\mathbf{E}_{i}\cdot\mathbf{E}_{j}\right)A^{i}A^{j}}=\sqrt{g_{ij}A^{i}A^{j}}=\sqrt{A_{j}A^{j}}=\sqrt{A^{i}A_{i}}\label{eqVectorMagnitudeCon}
\end{equation}
A similar expression can be obtained for the covariant form of the
vector, that is:
\begin{equation}
\left|\mathbf{A}\right|=\sqrt{\mathbf{A}\cdot\mathbf{A}}=\sqrt{\left(\mathbf{E}^{i}\cdot\mathbf{E}^{j}\right)A_{i}A_{j}}=\sqrt{g^{ij}A_{i}A_{j}}=\sqrt{A^{j}A_{j}}=\sqrt{A_{i}A^{i}}\label{eqVectorMagnitudeCov}
\end{equation}
The magnitude of a vector can also be obtained more directly from
the dot product of the covariant and contravariant forms of the vector:
\begin{equation}
\left|\mathbf{A}\right|=\sqrt{\mathbf{A}\cdot\mathbf{A}}=\sqrt{\left(\mathbf{E}^{i}\cdot\mathbf{E}_{j}\right)A_{i}A^{j}}=\sqrt{\delta_{j}^{i}A_{i}A^{j}}=\sqrt{A_{i}A^{i}}=\sqrt{A_{j}A^{j}}\label{eqVectorMagnitudeMix}
\end{equation}

\subsubsection{Angle Between Vectors}

$\bullet$ The angle $\theta$ between two contravariant or two covariant
vectors $\mathbf{A}$ and $\mathbf{B}$ is given respectively by:
\begin{equation}
\cos\theta=\frac{\mathbf{A}\cdot\mathbf{B}}{\left|\mathbf{A}\right|\left|\mathbf{B}\right|}=\frac{g_{ij}A^{i}B^{j}}{\sqrt{g_{kl}A^{k}A^{l}}\sqrt{g_{mn}B^{m}B^{n}}}=\frac{g^{ij}A_{i}B_{j}}{\sqrt{g^{kl}A_{k}A_{l}}\sqrt{g^{mn}B_{m}B_{n}}}
\end{equation}
For two vectors of opposite variance type we have:
\begin{equation}
\cos\theta=\frac{\mathbf{A}\cdot\mathbf{B}}{\left|\mathbf{A}\right|\left|\mathbf{B}\right|}=\frac{A^{i}B_{i}}{\sqrt{g_{kl}A^{k}A^{l}}\sqrt{g^{mn}B_{m}B_{n}}}=\frac{A_{i}B^{i}}{\sqrt{g^{kl}A_{k}A_{l}}\sqrt{g_{mn}B^{m}B^{n}}}
\end{equation}

\subsubsection{Length of Curve}

$\bullet$ In general curvilinear coordinates, the length of a $t$-parameterized
space curve $\mathbf{r}(t)$ defined by $u^{i}=u^{i}(t)$, which represents
the distance traversed along the curve on moving between its start
point $S$ and end point $E$, is given by:\footnote{Some authors add a sign indicator to ensure that the argument of the
square root is positive. However, as indicated in the Preface, such
a condition is assumed when needed since we deal with non-complex
values only.}
\begin{equation}
L=\int_{S}^{E}\sqrt{g_{ij}du^{i}du^{j}}=\int_{t_{1}}^{t_{2}}\sqrt{g_{ij}\frac{du^{i}}{dt}\frac{du^{j}}{dt}}\,\,dt
\end{equation}
where $t$ is a scalar variable parameter, and $t_{1}$ and $t_{2}$
are the values of $t$ corresponding to the start and end points respectively.

$\bullet$ The length of curve is used to define the geodesic which
is the path of the shortest distance connecting two points in a Riemannian
space. Although the geodesic is a straight line in a Euclidean space,
it is a generalized curved path in a general Riemannian space.

\pagebreak{}

\section{Covariant and Absolute Differentiation\label{secCovariantDifferentiation}}

$\bullet$ The focus of this section is the investigation of covariant
and absolute differentiation operations which are closely linked.
These operations represent generalization of tensor differentiation
in general curvilinear coordinate systems. Briefly, the differential
change of a tensor in general curvilinear coordinate systems is the
result of a change in the base vectors and a change in the tensor
components. Hence, covariant and absolute differentiation, in place
of the normal differentiation, are defined and employed to account
for both of these changes. Since Christoffel symbols are crucial in
the formulation and application of covariant and absolute differentiation,
the first subsection of the present section will be dedicated to these
symbols and their properties.

\subsection{Christoffel Symbols\label{subChristoffelSymbols}}

$\bullet$ We start by investigating the main properties of the Christoffel
symbols which play crucial roles in tensor calculus in general and
are needed for the subsequent development of the present and forthcoming
sections as well as the future notes.

$\bullet$ Christoffel symbols are classified as those of the first
kind and those of the second kind. These two kinds are linked through
the index raising and lowering operators. Both kinds of Christoffel
symbols are variable functions of coordinates in general.

$\bullet$ Christoffel symbols of the first and second kind are not
tensors in general although they are affine tensors of rank-3.

$\bullet$ As a consequence of the last point, if all the Christoffel
symbols of either kind vanished in a particular coordinate system
they will not necessarily vanish in other systems; for instance they
all vanish in Cartesian systems but not in cylindrical or spherical
systems, as has been established previously \cite{SochiTensorIntro2016}
and will be investigated further in the forthcoming points.

$\bullet$ Christoffel symbols of the first kind are given by:
\begin{equation}
\left[ij,l\right]=\frac{1}{2}\left(\partial_{j}g_{il}+\partial_{i}g_{jl}-\partial_{l}g_{ij}\right)
\end{equation}
where the indexed $g$ is the covariant form of the metric tensor.

$\bullet$ Christoffel symbols of the second kind are obtained by
raising the third index of the Christoffel symbols of the first kind,
that is:
\begin{equation}
\Gamma_{ij}^{k}=g^{kl}\left[ij,l\right]=\frac{g^{kl}}{2}\left(\partial_{j}g_{il}+\partial_{i}g_{jl}-\partial_{l}g_{ij}\right)\label{eqChris2nd}
\end{equation}
where the indexed $g$ is the metric tensor in its contravariant and
covariant forms with implied summation over $l$.

$\bullet$ Similarly, the Christoffel symbols of the first kind can
be obtained from the Christoffel symbols of the second kind by reversing
the above process through lowering the upper index, that is:
\begin{equation}
g_{km}\Gamma_{ij}^{k}=g_{km}g^{kl}\left[ij,l\right]=\delta_{m}^{l}\left[ij,l\right]=\left[ij,m\right]
\end{equation}

$\bullet$ For an $n$D space with $n$ covariant basis vectors ($\mathbf{E}_{1},\mathbf{E}_{2},\ldots,\mathbf{E}_{n}$)
spanning the space, the derivative $\partial_{j}\mathbf{E}_{i}$ for
any given $i$ is a vector within the space and hence it is in general
a linear combination of all the basis vectors. The Christoffel symbols
of the second kind are the components of this linear combination,
that is:
\begin{equation}
\partial_{j}\mathbf{E}_{i}=\Gamma_{ij}^{k}\mathbf{E}_{k}\label{eqPartialBasisCovariant}
\end{equation}
Similarly, for the contravariant basis vectors we have:
\begin{equation}
\partial_{j}\mathbf{E}^{i}=-\Gamma_{kj}^{i}\mathbf{E}^{k}\label{eqPartialBasisContravariant}
\end{equation}

$\bullet$ By inner product multiplication of the previous relations
with the basis vectors we obtain:
\begin{equation}
\mathbf{E}^{i}\cdot\partial_{k}\mathbf{E}_{j}=\Gamma_{jk}^{i}\,\,\,\,\,\,\,\,\,\,\,\,\,\,\,\,\,\&\,\,\,\,\,\,\,\,\,\,\,\,\,\,\,\,\mathbf{E}_{i}\cdot\partial_{k}\mathbf{E}^{j}=-\Gamma_{ik}^{j}
\end{equation}
Similarly:
\begin{equation}
\mathbf{E}_{k}\cdot\partial_{j}\mathbf{E}_{i}=g_{mk}\mathbf{E}^{m}\cdot\partial_{j}\mathbf{E}_{i}=g_{mk}\Gamma_{ij}^{m}=\left[ij,k\right]
\end{equation}

$\bullet$ Christoffel symbols of the first and second kind are symmetric
in their paired indices, that is:
\begin{equation}
\left[ij,k\right]=\left[ji,k\right]\,\,\,\,\,\,\,\,\,\,\,\,\,\,\,\,\&\,\,\,\,\,\,\,\,\,\,\,\,\,\,\,\,\,\,\Gamma_{ij}^{k}=\Gamma_{ji}^{k}
\end{equation}

$\bullet$ The partial derivative of the components of the covariant
metric tensor and the Christoffel symbols of the first kind satisfy
the following identity, which is essentially based on the forthcoming
Ricci Theorem:
\begin{equation}
\partial_{j}g_{il}=\left[ij,l\right]+\left[jl,i\right]\label{eqChristoffelFirstKindRelation}
\end{equation}
This relation can also be written in terms of the Christoffel symbols
of the second kind using the index shifting operator:
\begin{equation}
\partial_{j}g_{il}=g_{kl}\Gamma_{ij}^{k}+g_{ki}\Gamma_{jl}^{k}
\end{equation}
A related formula for the partial derivative of the components of
the contravariant metric tensor, which can be obtained by partial
differentiation of the relation $g_{im}g^{mj}=\delta_{i}^{j}$ with
respect to the $k^{th}$ coordinate, is given by:
\begin{equation}
g_{im}\partial_{k}g^{mj}=-g^{mj}\partial_{k}g_{im}
\end{equation}

$\bullet$ Christoffel symbols of the second kind with two identical
indices of opposite variance type satisfy the following relations:
\begin{equation}
\Gamma_{ji}^{j}=\Gamma_{ij}^{j}=\frac{1}{2g}\partial_{i}g=\frac{1}{2}\partial_{i}\left(\ln\,g\right)=\partial_{i}\left(\ln\sqrt{g}\right)=\frac{1}{\sqrt{g}}\partial_{i}\sqrt{g}\label{eqGammaaia}
\end{equation}
where the main relation can be derived as follow:
\begin{equation}
\begin{aligned}\Gamma_{ij}^{j} & =\frac{g^{jl}}{2}\left(\partial_{j}g_{il}+\partial_{i}g_{jl}-\partial_{l}g_{ij}\right) & \,\,\,\,\,\, & \text{(Eq. \ref{eqChris2nd} with \ensuremath{k=j})}\\
 & =\frac{g^{jl}}{2}\left(\partial_{l}g_{ij}+\partial_{i}g_{jl}-\partial_{l}g_{ij}\right) &  & \text{(relabeling dummy \ensuremath{j,l} in 1\ensuremath{^{st}} term \& \ensuremath{g^{jl}=g^{lj}})}\\
 & =\frac{1}{2}g^{jl}\partial_{i}g_{jl}\\
 & =\frac{1}{2g}gg^{jl}\partial_{i}g_{jl}\\
 & =\frac{1}{2g}\partial_{i}g &  & \text{(derivative of determinant)}
\end{aligned}
\end{equation}

$\bullet$ In orthogonal coordinate systems, the Christoffel symbols
of the first kind are given by:
\begin{equation}
\begin{aligned}\left[ij,i\right] & =\left[ji,i\right]=\frac{1}{2}\partial_{j}g_{ii} & \,\,\,\,\,\,\,\,\,\,\,\,\,\,\,\,\, & \text{(no sum on \ensuremath{i})}\\
\left[ii,j\right] & =-\frac{1}{2}\partial_{j}g_{ii} &  & \text{(\ensuremath{i\ne j}, no sum on \ensuremath{i})}\\
\left[ij,k\right] & =0 &  & (\ensuremath{i\ne j\ne k})
\end{aligned}
\end{equation}
The first relation is a special case of Eq. \ref{eqChristoffelFirstKindRelation}
with $l=i$ taking into account that the Christoffel symbols are symmetric
in their paired indices; moreover, the relation includes the case
of $i=j$, i.e. when all the three indices are identical.

$\bullet$ In orthogonal coordinate systems, the Christoffel symbols
of the second kind are given by:
\begin{equation}
\Gamma_{jk}^{i}=\frac{\left[jk,i\right]}{g_{ii}}=g^{ii}\left[jk,i\right]\,\,\,\,\,\,\,\,\,\,\,\,\,\,\,\,\,\,\,\,\text{(no sum on \ensuremath{i})}
\end{equation}
and hence from the results of the previous point we have:
\begin{equation}
\begin{aligned}\Gamma_{ij}^{i} & =\Gamma_{ji}^{i}=\frac{1}{2g_{ii}}\partial_{j}g_{ii}=\frac{g^{ii}}{2}\partial_{j}g_{ii}=\frac{1}{2}\partial_{j}\ln g_{ii} & \,\,\,\,\,\,\,\,\,\, & \text{(no sum on \ensuremath{i})}\\
\Gamma_{jj}^{i} & =-\frac{1}{2g_{ii}}\partial_{i}g_{jj}=-\frac{g^{ii}}{2}\partial_{i}g_{jj} &  & \text{(no sum on \ensuremath{i} or \ensuremath{j}, and \ensuremath{i\ne j})}\\
\Gamma_{jk}^{i} & =0 &  & (\ensuremath{i\ne j\ne k})
\end{aligned}
\end{equation}
As for the first kind in the last point, the first relation includes
the case of $i=j$, i.e. when all the three indices are identical.

$\bullet$ In general orthogonal coordinate systems in a 3D space,
the Christoffel symbols of the first kind vanish when the indices
are all different, as shown earlier; moreover, the non-identically
vanishing symbols of the first kind are given by:
\begin{eqnarray}
[11,1]=+h_{1}h_{1,1} & \,\,\,\,\,\,\,\,\,\,\,\,\,\,\,[11,2]=-h_{1}h_{1,2} & \,\,\,\,\,\,\,\,\,\,\,\,\,\,\,[11,3]=-h_{1}h_{1,3}\nonumber \\{}
[12,1]=+h_{1}h_{1,2} & \,\,\,\,\,\,\,\,\,\,\,\,\,\,\,[12,2]=+h_{2}h_{2,1} & \,\,\,\,\,\,\,\,\,\,\,\,\,\,\,[13,1]=+h_{1}h_{1,3}\nonumber \\{}
[13,3]=+h_{3}h_{3,1} & \,\,\,\,\,\,\,\,\,\,\,\,\,\,\,[22,1]=-h_{2}h_{2,1} & \,\,\,\,\,\,\,\,\,\,\,\,\,\,\,[22,2]=+h_{2}h_{2,2}\\{}
[22,3]=-h_{2}h_{2,3} & \,\,\,\,\,\,\,\,\,\,\,\,\,\,\,[23,2]=+h_{2}h_{2,3} & \,\,\,\,\,\,\,\,\,\,\,\,\,\,\,[23,3]=+h_{3}h_{3,2}\nonumber \\{}
[33,1]=-h_{3}h_{3,1} & \,\,\,\,\,\,\,\,\,\,\,\,\,\,\,[33,2]=-h_{3}h_{3,2} & \,\,\,\,\,\,\,\,\,\,\,\,\,\,\,[33,3]=+h_{3}h_{3,3}\nonumber
\end{eqnarray}
where ($1,2,3$) stand for ($u^{1},$$u^{2},$$u^{3}$) respectively,
$h_{1},h_{2},h_{3}$ are the scale factors as defined previously,
and the comma indicates, as always, partial derivative; for example
in cylindrical coordinates given by ($\rho,\phi,z$), $h_{2,1}$ means
the partial derivative of $h_{2}$ with respect to the first coordinate
and hence $h_{2,1}=\partial_{\rho}\rho=1$ since $h_{2}=\rho$ and
the first coordinate is $\rho$ (refer to Table \ref{tabScaleFactors}).
Because the Christoffel symbols of the first kind are symmetric in
their first two indices, the $\left[21,1\right]$ symbol for instance
can be obtained from the value of the $\left[12,1\right]$ symbol.

$\bullet$ In general orthogonal coordinate systems in a 3D space,
the Christoffel symbols of the second kind vanish when the indices
are all different, as shown earlier; moreover, the non-identically
vanishing symbols of the second kind are given by:
\begin{eqnarray}
\Gamma_{11}^{1}=+\frac{h_{1,1}}{h_{1}}\,\,\,\,\,\, & \,\,\,\,\,\,\,\,\,\,\,\,\,\,\,\Gamma_{11}^{2}=-\frac{h_{1}h_{1,2}}{\left(h_{2}\right)^{2}} & \,\,\,\,\,\,\,\,\,\,\,\,\,\,\,\Gamma_{11}^{3}=-\frac{h_{1}h_{1,3}}{\left(h_{3}\right)^{2}}\nonumber \\
\Gamma_{12}^{1}=+\frac{h_{1,2}}{h_{1}}\,\,\,\,\,\, & \,\,\,\,\,\,\,\,\,\,\,\,\,\,\,\Gamma_{12}^{2}=+\frac{h_{2,1}}{h_{2}}\,\,\,\,\,\, & \,\,\,\,\,\,\,\,\,\,\,\,\,\,\,\Gamma_{13}^{1}=+\frac{h_{1,3}}{h_{1}}\nonumber \\
\Gamma_{13}^{3}=+\frac{h_{3,1}}{h_{3}}\,\,\,\,\,\, & \,\,\,\,\,\,\,\,\,\,\,\,\,\,\,\Gamma_{22}^{1}=-\frac{h_{2}h_{2,1}}{\left(h_{1}\right)^{2}} & \,\,\,\,\,\,\,\,\,\,\,\,\,\,\,\Gamma_{22}^{2}=+\frac{h_{2,2}}{h_{2}}\\
\Gamma_{22}^{3}=-\frac{h_{2}h_{2,3}}{\left(h_{3}\right)^{2}} & \,\,\,\,\,\,\,\,\,\,\,\,\,\,\,\Gamma_{23}^{2}=+\frac{h_{2,3}}{h_{2}}\,\,\,\,\,\, & \,\,\,\,\,\,\,\,\,\,\,\,\,\,\,\Gamma_{23}^{3}=+\frac{h_{3,2}}{h_{3}}\nonumber \\
\Gamma_{33}^{1}=-\frac{h_{3}h_{3,1}}{\left(h_{1}\right)^{2}} & \,\,\,\,\,\,\,\,\,\,\,\,\,\,\,\Gamma_{33}^{2}=-\frac{h_{3}h_{3,2}}{\left(h_{2}\right)^{2}} & \,\,\,\,\,\,\,\,\,\,\,\,\,\,\,\Gamma_{33}^{3}=+\frac{h_{3,3}}{h_{3}}\nonumber
\end{eqnarray}
where ($1,2,3$) stand for ($u^{1},$$u^{2},$$u^{3}$) respectively.
Again, since the Christoffel symbols of the second kind are symmetric
in their lower indices, the missing non-vanishing entries can be obtained
from the given entries by permuting the lower indices.

$\bullet$ In Cartesian coordinate systems ($x,y,z$), all the Christoffel
symbols of the first and second kind are identically zero.

$\bullet$ In cylindrical coordinate systems ($\rho,\phi,z$), the
non-zero Christoffel symbols of the first kind are:
\begin{eqnarray}
\left[22,1\right] & = & -\rho\\
\left[12,2\right] & = & \begin{aligned}\left[21,2\right]\,\, & \begin{aligned}= & \,\,\,\,\rho\end{aligned}
\end{aligned}
\nonumber
\end{eqnarray}
where ($1,2,3$) stand for ($\rho,\phi,z$) respectively.

$\bullet$ In cylindrical coordinate systems ($\rho,\phi,z$), the
non-zero Christoffel symbols of the second kind are:
\begin{eqnarray}
\Gamma_{22}^{1} & = & -\rho\\
\Gamma_{12}^{2} & = & \begin{aligned}\Gamma_{21}^{2} & \begin{aligned}\,\,=\,\,\,\, & \frac{1}{\rho}\end{aligned}
\end{aligned}
\nonumber
\end{eqnarray}
where ($1,2,3$) stand for ($\rho,\phi,z$) respectively.

$\bullet$ In spherical coordinate systems ($r,\theta,\phi$), the
non-zero Christoffel symbols of the first kind are:
\begin{eqnarray}
\left[22,1\right] & = & -r\\
\left[33,1\right] & = & -r\sin^{2}\theta\nonumber \\
\left[12,2\right] & = & \begin{aligned}\left[21,2\right] & \begin{aligned}\,=\,\,\,\,\,\, & r\end{aligned}
\end{aligned}
\nonumber \\
\left[33,2\right] & = & -r^{2}\sin\theta\cos\theta\nonumber \\
\left[13,3\right] & = & \begin{aligned}\begin{aligned}\left[31,3\right] & =\end{aligned}
 & r\sin^{2}\theta\end{aligned}
\nonumber \\
\left[23,3\right] & = & \begin{aligned}\begin{aligned}\left[32,3\right] & =\end{aligned}
 & r^{2}\sin\theta\cos\theta\end{aligned}
\nonumber
\end{eqnarray}
where ($1,2,3$) stand for ($r,\theta,\phi$) respectively.

$\bullet$ In spherical coordinate systems ($r,\theta,\phi$), the
non-zero Christoffel symbols of the second kind are:
\begin{eqnarray}
\Gamma_{22}^{1} & = & -r\\
\Gamma_{33}^{1} & = & -r\sin^{2}\theta\nonumber \\
\Gamma_{12}^{2} & = & \begin{aligned}\Gamma_{21}^{2} & \begin{aligned}\,=\,\,\,\,\,\, & \frac{1}{r}\end{aligned}
\end{aligned}
\nonumber \\
\Gamma_{33}^{2} & = & -\sin\theta\cos\theta\nonumber \\
\Gamma_{13}^{3} & = & \begin{aligned}\begin{aligned}\Gamma_{31}^{3} & =\end{aligned}
 & \frac{1}{r}\end{aligned}
\nonumber \\
\Gamma_{23}^{3} & = & \begin{aligned}\begin{aligned}\Gamma_{32}^{3} & =\end{aligned}
 & \cot\theta\end{aligned}
\nonumber
\end{eqnarray}
where ($1,2,3$) stand for ($r,\theta,\phi$) respectively.

$\bullet$ Because there is an element of arbitrariness in the choice
of the coordinates order and hence their indices, the Christoffel
symbols may be given in terms of coordinate symbols rather than their
indices to be more explicit and avoid ambiguity; for instance in the
above examples of cylindrical and spherical coordinate systems we
have: $\left[22,1\right]\equiv\left[\phi\phi,\rho\right]$, $\Gamma_{12}^{2}\equiv\Gamma_{\rho\phi}^{\phi}$
for cylindrical, and $\left[22,1\right]\equiv\left[\theta\theta,r\right]$,
$\Gamma_{13}^{3}\equiv\Gamma_{r\phi}^{\phi}$ for spherical.

$\bullet$ The Christoffel symbols may be subscripted by the symbol
of the metric tensor for the given space to reveal the metric which
they are based upon.

$\bullet$ In any coordinate system, all the Christoffel symbols of
the first and second kind vanish identically \textit{iff} all the
components of the metric tensor in the given coordinate system are
constants.

$\bullet$ In affine coordinates, all the components of the metric
tensor are constants and hence all the Christoffel symbols of both
kinds vanish identically.

$\bullet$ The number of independent Christoffel symbols of each kind
(first and second) in general curvilinear coordinates is given by:
\begin{equation}
N_{\mathrm{CI}}=\frac{n^{2}\left(n+1\right)}{2}
\end{equation}
where $n$ is the space dimension. The reason is that, due to the
symmetry of the metric tensor there are $\frac{n\left(n+1\right)}{2}$
independent metric components, $g_{ij}$, and for each independent
component there are $n$ distinct Christoffel symbols.

$\bullet$ The following relations are useful in the manipulation
of tensor expressions involving Christoffel symbols:\footnote{The first relation is a special case of the relation: $A_{\,\,;j}^{ij}=\frac{1}{\sqrt{g}}\partial_{j}\left(\sqrt{g}A^{ij}\right)+A^{kl}\Gamma_{kl}^{i}$
noting that the covariant derivative of the metric tensor is identically
zero according to the Ricci Theorem.}
\begin{equation}
\frac{1}{\sqrt{g}}\partial_{j}\left(\sqrt{g}g^{ij}\right)+g^{kl}\Gamma_{kl}^{i}=0
\end{equation}
\begin{equation}
\partial_{j}\left[ik,l\right]=g_{la}\partial_{j}\Gamma_{ik}^{a}+\Gamma_{ik}^{a}\left[lj,a\right]+\Gamma_{ik}^{a}\left[aj,l\right]
\end{equation}

\subsection{Covariant Derivative\label{subCovariantDerivative}}

$\bullet$ The basis vectors in general curvilinear coordinate systems
undergo changes in magnitude and direction as they move around in
their own space, and hence they are functions of position. These changes
should be accounted for when calculating the derivatives of tensors
in such general systems. Therefore, terms based on using Christoffel
symbols are added to the ordinary derivative terms to correct for
these changes and this more comprehensive form of derivative is called
the covariant derivative.

$\bullet$ Since in rectilinear coordinate systems the basis vectors
are constants, the Christoffel symbol terms vanish identically and
hence the covariant derivative reduces to the ordinary derivative,
but in the other coordinate systems these terms are present in general.

$\bullet$ As a consequence of the last point, the ordinary derivative
of a non-scalar tensor is a tensor \textit{iff} the coordinate transformations
are linear.

$\bullet$ It has been stated that the ``covariant'' label is an
indication that the differentiation operator, $\nabla_{;i}$, is in
the covariant position. However, it may also be true that ``covariant''
means ``invariant'' as pointed out earlier in the previous set of
notes.

$\bullet$ Contravariant differentiation ($\nabla^{;j}$) can also
be defined for covariant and contravariant tensors by raising the
differentiation index using the index raising operator, e.g.
\begin{equation}
A_{i}^{\,\,;j}=g^{jk}A_{i;k}\,\,\,\,\,\,\,\,\,\,\,\,\,\,\,\,\,\&\,\,\,\,\,\,\,\,\,\,\,\,\,\,\,\,\,A^{i;j}=g^{jk}A_{\,\,\,;k}^{i}
\end{equation}
However, practically such operations are rarely used.\footnote{An example of contravariant differentiation is in the definition of
the Laplacian in general curvilinear coordinates (refer to $\S$ \ref{subLaplacianGneralCurvilinear}).}

$\bullet$ As an example of how to obtain the covariant derivative
of a tensor, let have a vector represented by contravariant components:
$\mathbf{A}=A^{i}\mathbf{E}_{i}$ in general coordinates. We differentiate
this vector following the normal rules of differentiation and taking
account of the fact that the basis vectors in general coordinates
are differentiable functions of position and hence they, unlike their
Cartesian counterparts, are subject to differentiation using the product
rule, that is:
\begin{equation}
\begin{aligned}\partial_{j}\mathbf{A} & =\mathbf{E}_{i}\partial_{j}A^{i}+A^{i}\partial_{j}\mathbf{E}_{i} & \,\,\,\,\,\,\,\,\,\,\,\,\, & \text{(product rule)}\\
 & =\mathbf{E}_{i}\partial_{j}A^{i}+A^{i}\Gamma_{\,\,\,ij}^{k}\mathbf{E}_{k} &  & \text{(Eq. \ref{eqPartialBasisCovariant})}\\
 & =\mathbf{E}_{i}\partial_{j}A^{i}+A^{k}\Gamma_{\,\,\,kj}^{i}\mathbf{E}_{i} &  & \text{(relabeling dummy indices \ensuremath{i} \& \ensuremath{k})}\\
 & =\left(\partial_{j}A^{i}+A^{k}\Gamma_{\,\,\,kj}^{i}\right)\mathbf{E}_{i}\\
 & =A_{\,\,;j}^{i}\mathbf{E}_{i}
\end{aligned}
\label{eqCoDifDev}
\end{equation}
where $A_{\,\,;j}^{i}$, which is a rank-2 mixed tensor, is labeled
the ``covariant derivative'' of $A^{i}$.

Similarly, for a vector represented by covariant components: $\mathbf{A}=A_{i}\mathbf{E}^{i}$
in general curvilinear coordinates we have:
\begin{equation}
\partial_{j}\mathbf{A}=A_{i;j}\mathbf{E}^{i}\label{eqCoDifCo}
\end{equation}

$\bullet$ Following the method and techniques outlined in the previous
point, to obtain the covariant derivative of a tensor in general,
we start with an ordinary partial derivative term of the given tensor.
Then for each tensor index an extra Christoffel symbol term is added,
positive for contravariant indices and negative for covariant indices,
where the differentiation index is one of the lower indices in the
Christoffel symbol. Hence, for a general differentiable rank-$n$
tensor $\mathbf{A}$ the covariant derivative is given by:
\begin{eqnarray}
A_{lm\ldots p;q}^{ij\ldots k} & = & \partial_{q}A_{lm\ldots p}^{ij\ldots k}+\Gamma_{aq}^{i}A_{lm\ldots p}^{aj\ldots k}+\Gamma_{aq}^{j}A_{lm\ldots p}^{ia\ldots k}+\cdots+\Gamma_{aq}^{k}A_{lm\ldots p}^{ij\ldots a}\label{eqCovariantDerivative}\\
 &  & \,\,\,\,\,\,\,\,\,\,\,\,\,\,\,\,\,\,\,\,\,\,-\Gamma_{lq}^{a}A_{am\ldots p}^{ij\ldots k}-\Gamma_{mq}^{a}A_{la\ldots p}^{ij\ldots k}-\cdots-\Gamma_{pq}^{a}A_{lm\ldots a}^{ij\ldots k}\nonumber
\end{eqnarray}

$\bullet$ Practically, there is only one possibility for the arrangement
of the indices in the Christoffel symbol terms if the following rules
are observed:

(A) the second subscript index of the Christoffel symbol is the differentiation
index,

(B) the concerned tensor index in the Christoffel symbol term is contracted
with one of the indices of the Christoffel symbol and hence they are
opposite in their lower/upper position,

(C) the contracted index is transferred from the tensor to the Christoffel
symbol keeping its lower/upper position, and

(D) all the other indices of the tensor keep their names and position.

$\bullet$ The ordinary partial derivative term in the above covariant
derivative expression (Eq. \ref{eqCovariantDerivative}) represents
the rate of change of the tensor components with change of position
as a result of moving along the coordinate curve of the differentiated
index, while the Christoffel symbol terms represent the change experienced
by the local basis vectors as a result of the same movement. This
can be seen from the development of Eq. \ref{eqCoDifDev}.

$\bullet$ From the above discussion it is obvious that to obtain
the covariant derivative, the Christoffel symbols are required and
these symbols are dependent on the metric tensor; hence the covariant
derivative is dependent on having the space metric.

$\bullet$ In all coordinate systems, the covariant derivative of
a differentiable scalar function of position, $f$, is the same as
the ordinary partial derivative, that is:
\begin{equation}
f_{;i}=f_{,i}=\partial_{i}f
\end{equation}
This is justified by the fact that the covariant derivative is different
from the ordinary partial derivative because the basis vectors in
general coordinate systems are dependent on their spatial position,
and since a scalar is independent of the basis vectors the covariant
derivative and partial derivative are identical. This can also be
concluded from the covariant derivative rule as stated in the previous
points and formulated in Eq. \ref{eqCovariantDerivative}.

$\bullet$ Several rules of normal differentiation are naturally extended
to covariant differentiation. For example, covariant differentiation
is a linear operation with respect to algebraic sums of tensor terms
and hence the covariant derivative of a sum is the sum of the covariant
derivatives of the terms:
\begin{equation}
\left(a\mathbf{A}\pm b\mathbf{B}\right)_{;i}=a\left(\mathbf{A}\right)_{;i}\pm b\left(\mathbf{B}\right)_{;i}
\end{equation}
where $a$ and $b$ are scalar constants and $\mathbf{A}$ and $\mathbf{B}$
are differentiable tensors. The product rule of ordinary differentiation
also applies to covariant differentiation of inner and outer products
of tensors:
\begin{equation}
\left(\mathbf{A}\circ\mathbf{B}\right)_{;i}=\left(\mathbf{A}\right)_{;i}\circ\mathbf{B}+\mathbf{A}\circ\left(\mathbf{B}\right)_{;i}
\end{equation}
where the symbol $\circ$ denotes an inner or outer product operator.

$\bullet$ According to the ``Ricci Theorem'', the covariant derivative
of the covariant and contravariant metric tensor is zero. This has
nothing to do with the metric tensor being a constant function of
coordinates, which is true only for the rectilinear systems, but this
arises from the fact that the covariant derivative quantifies the
change with position of the basis vectors in magnitude and direction
as well as the change in components, and these contributions in the
case of the metric tensor cancel each other resulting in a total null
effect. As a result, the metric tensor behaves as a constant with
respect to the covariant derivative operation:
\begin{equation}
g_{ij;k}=0\,\,\,\,\,\,\,\,\,\,\,\,\,\,\,\,\,\&\,\,\,\,\,\,\,\,\,\,\,\,\,\,\,\,\,g_{;k}^{ij}=0
\end{equation}
for all values of the indices, and hence the covariant derivative
operator bypasses the metric tensor:
\begin{equation}
\left(\mathbf{g}\circ\mathbf{A}\right)_{;k}=\mathbf{g}\circ\mathbf{\left(A\right)}_{;k}\label{eqCovDerivMetric}
\end{equation}
where $\mathbf{A}$ is a general tensor, $\mathbf{g}$ is the metric
tensor in its covariant or contravariant form and $\circ$ denotes
an inner or outer tensor product.\footnote{Although the metric tensor is normally used in inner product operations
for raising and lowering of indices, the possibility of its involvement
in outer product operations should not be ruled out.}

$\bullet$ As a result of the Ricci Theorem, the covariant derivative
operator and the index shifting operator are commutative, e.g.
\begin{equation}
\left(g_{ik}A^{k}\right)_{;j}=g_{ik}A_{;j}^{k}=A_{i;j}=\left(A_{i}\right)_{;j}=\left(g_{ik}A^{k}\right)_{;j}
\end{equation}
\begin{equation}
\left(g_{mi}g_{nj}A^{mn}\right)_{;k}=g_{mi}g_{nj}A_{\,;k}^{mn}=A_{ij;k}=\left(A_{ij}\right)_{;k}=\left(g_{mi}g_{nj}A^{mn}\right)_{;k}
\end{equation}

$\bullet$ Like the metric tensor, the Kronecker delta is constant
with regard to the covariant differentiation and hence the covariant
derivative of the Kronecker delta is identically zero:
\begin{equation}
\delta_{j;k}^{i}=\partial_{k}\delta_{j}^{i}+\delta_{j}^{a}\Gamma_{ak}^{i}-\delta_{a}^{i}\Gamma_{jk}^{a}=0+\Gamma_{jk}^{i}-\Gamma_{jk}^{i}=0
\end{equation}
The rule of the Kronecker delta may be regarded as an instance of
the rule of the metric tensor, as stated by the Ricci Theorem, since
the Kronecker delta is a metric tensor. Likewise, the covariant differentiation
operator bypasses the Kronecker delta which is involved in inner and
outer tensor products:\footnote{Like the metric tensor, the Kronecker delta is normally used in inner
product operations for replacement of indices; however the possibility
of its involvement in outer product operations should not be ruled
out.}
\begin{equation}
\left(\boldsymbol{\delta}\circ\mathbf{A}\right)_{;k}=\boldsymbol{\delta}\circ\mathbf{\left(A\right)}_{;k}
\end{equation}

$\bullet$ Like the ordinary Kronecker delta, the covariant derivative
of the generalized Kronecker delta is identically zero.

$\bullet$ For a differentiable function $f(x,y)$ of class $C^{2}$
(i.e. all the second order partial derivatives of the function do
exist and are continuous), the mixed partial derivatives are equal,
that is:
\begin{equation}
\partial_{x}\partial_{y}f=\partial_{y}\partial_{x}f
\end{equation}
However, even if the components of a tensor satisfy this condition
(i.e. being of class $C^{2}$), this is not sufficient for the equality
of the mixed covariant derivatives. What is required for the mixed
covariant derivatives to be equal is the vanishing of the Riemann
Tensor (see $\S$ \ref{subRiemannTensor}).

$\bullet$ Higher order covariant derivatives are similarly defined
as derivatives of derivatives by successive repetition of the process
of covariant differentiation; however the order of differentiation
should be respected as stated in the previous point. For example,
the second order mixed $jk$ covariant derivative of a contravariant
vector $\mathbf{A}$ is given by:
\begin{equation}
A_{;jk}^{i}=\partial_{k}\text{\ensuremath{\partial}}_{j}A^{i}+\Gamma_{ka}^{i}\text{\ensuremath{\partial}}_{j}A^{a}-\Gamma_{jk}^{a}\text{\ensuremath{\partial}}_{a}A^{i}+\text{\ensuremath{\Gamma}}_{ja}^{i}\text{\ensuremath{\partial}}_{k}A^{a}+A^{a}\left(\text{\ensuremath{\partial}}_{j}\Gamma_{ka}^{i}-\text{\ensuremath{\Gamma}}_{jk}^{b}\text{\ensuremath{\Gamma}}_{ba}^{i}+\text{\ensuremath{\Gamma}}_{jb}^{i}\text{\ensuremath{\Gamma}}_{ka}^{b}\right)
\end{equation}
while the second order mixed $jk$ covariant derivative of a covariant
vector $\mathbf{A}$ is given by:
\begin{equation}
A_{i;jk}=\partial_{k}\text{\ensuremath{\partial}}_{j}A_{i}-\Gamma_{ij}^{a}\text{\ensuremath{\partial}}_{k}A_{a}-\Gamma_{ik}^{a}\text{\ensuremath{\partial}}_{j}A_{a}-\text{\ensuremath{\Gamma}}_{jk}^{a}\text{\ensuremath{\partial}}_{a}A_{i}-A_{a}\left(\text{\ensuremath{\partial}}_{k}\Gamma_{ij}^{a}-\text{\ensuremath{\Gamma}}_{ib}^{a}\text{\ensuremath{\Gamma}}_{jk}^{b}-\text{\ensuremath{\Gamma}}_{ik}^{b}\text{\ensuremath{\Gamma}}_{bj}^{a}\right)
\end{equation}

$\bullet$ The second order mixed $kj$ covariant derivative of contravariant
and covariant vectors can be obtained from the equations in the last
point by interchanging the $j$ and $k$ indices and hence the inequality
of the $jk$ and $kj$ mixed derivatives in general can be verified
(refer to $\S$ \ref{subRiemannTensor}).

$\bullet$ The covariant derivative of a tensor is a tensor whose
covariant rank is higher than the covariant rank of the original tensor
by one. Hence, the covariant derivative of a rank-$n$ tensor of type
($r,s$) is a rank-($n+1$) tensor of type ($r,s+1$).

$\bullet$ Covariant differentiation and contraction of index operations
commute with each other, e.g.
\begin{equation}
\left(A_{k}^{ij}\right)_{;l}\delta_{j}^{k}=\left(A_{k}^{ij}\delta_{j}^{k}\right)_{;l}=\left(A_{k}^{ik}\right)_{;l}
\end{equation}

$\bullet$ Since the Christoffel symbols vanish when the components
of the metric tensor in a given coordinate system are constants, the
covariant derivative is reduced to the ordinary derivative in such
systems. This is particularly true in Euclidean spaces coordinated
by rectilinear systems.

$\bullet$ The covariant derivatives of relative tensors, which are
also relative tensors of the same weight as the original tensors,
are obtained by adding a weight term to the normal formulae of covariant
derivative. Hence, the covariant derivative of a relative scalar with
weight $w$ is given by:
\begin{equation}
f_{;i}=f_{,i}-wf\Gamma_{ji}^{j}
\end{equation}
while the covariant derivative of relative tensors of higher ranks
with weight $w$ is obtained by adding the following term to the right
hand side of Eq. \ref{eqCovariantDerivative}:
\begin{equation}
-wA_{lm\ldots p}^{ij\ldots k}\Gamma_{aq}^{a}
\end{equation}

$\bullet$ Unlike ordinary differentiation, the covariant derivative
of a non-scalar tensor with constant components is not zero in general
due to the presence of the Christoffel symbols in the definition of
the covariant derivative, as given by Eq. \ref{eqCovariantDerivative}.

$\bullet$ In rectilinear coordinates, the Christoffel symbols are
identically zero because the basis vectors are constants, and hence
the covariant derivative is the same as the normal partial derivative
for all tensor ranks. As a result, when the components of the metric
tensor, $g_{ij}$, are constants as in the case of rectangular coordinate
systems, the covariant derivative becomes the ordinary partial derivative.

$\bullet$ For a differentiable covariant vector $\mathbf{A}$ which
is a gradient of a scalar we have:
\begin{equation}
A_{i;j}=A_{j;i}
\end{equation}

$\bullet$ The covariant derivative of the basis vectors of the covariant
and contravariant types is identically zero:
\begin{equation}
\begin{aligned}\mathbf{E}_{i;j} & =\partial_{j}\mathbf{E}_{i}-\Gamma_{ij}^{k}\mathbf{E}_{k}=\Gamma_{ij}^{k}\mathbf{E}_{k}-\Gamma_{ij}^{k}\mathbf{E}_{k}=\mathbf{0}\\
\mathbf{E}_{;j}^{i} & =\partial_{j}\mathbf{E}^{i}+\Gamma_{kj}^{i}\mathbf{E}^{k}=-\Gamma_{kj}^{i}\mathbf{E}^{k}+\Gamma_{kj}^{i}\mathbf{E}^{k}=\mathbf{0}
\end{aligned}
\end{equation}

\subsection{Absolute Derivative\label{subAbsoluteDerivative}}

$\bullet$ The absolute derivative of a tensor along a $t$-parameterized
curve in an $n$D space with respect to the parameter $t$ is the
inner product of the covariant derivative of the tensor and the tangent
vector to the curve. In brief, the absolute derivative is a covariant
derivative of a tensor along a curve.

$\bullet$ For a tensor $A^{i}$, the inner product of $A_{;j}^{i}$,
which is a tensor, with another tensor is a tensor. Now, if the other
tensor is $\frac{du^{i}}{dt}$, which is the tangent vector to a $t$-parameterized
curve $L$ given by the equations $u^{i}=u^{i}(t)$, then the inner
product:
\begin{equation}
A_{\,;r}^{i}\frac{du^{r}}{dt}\label{eqAbsoluteDerivative}
\end{equation}
is a tensor of the same rank and type as the tensor $A^{i}$. The
tensor given by the expression \ref{eqAbsoluteDerivative} is called
the ``absolute'' or ``intrinsic'' or ``absolute covariant''
derivative of the tensor $A^{i}$ along the curve $L$ and is symbolized
by:
\begin{equation}
\frac{\delta A^{i}}{\delta t}
\end{equation}

$\bullet$ For a differentiable scalar $f$, the absolute derivative,
like the covariant derivative, is the same as the ordinary derivative,
that is:
\begin{equation}
\frac{\delta f}{\delta t}=\frac{df}{dt}
\end{equation}

$\bullet$ The absolute derivative of a differentiable contravariant
vector $A^{k}$ with respect to the parameter $t$ is given by:
\begin{equation}
\frac{\delta A^{k}}{\delta t}=\frac{dA^{k}}{dt}+\Gamma_{ij}^{k}A^{i}\frac{du^{j}}{dt}
\end{equation}
Similarly for a differentiable covariant vector $A_{k}$ we have:
\begin{equation}
\frac{\delta A_{k}}{\delta t}=\frac{dA_{k}}{dt}-\Gamma_{kj}^{i}A_{i}\frac{du^{j}}{dt}
\end{equation}

$\bullet$ Absolute differentiation can be easily extended to higher
rank ($>1$) differentiable tensors of type ($m,n$) along parameterized
curves. For instance, the absolute derivative of a mixed tensor of
type ($1,2$) $A_{jk}^{i}$ along a $t$-parameterized curve $L$
is given by:
\begin{equation}
\frac{\delta A_{jk}^{i}}{\delta t}=A_{jk;b}^{i}\frac{du^{b}}{dt}=\frac{dA_{jk}^{i}}{dt}+\Gamma_{ab}^{i}A_{jk}^{a}\frac{du^{b}}{dt}-\Gamma_{jb}^{a}A_{ak}^{i}\frac{du^{b}}{dt}-\Gamma_{bk}^{a}A_{ja}^{i}\frac{du^{b}}{dt}
\end{equation}

$\bullet$ As the absolute derivative is given generically by:
\begin{equation}
\frac{\delta\mathbf{A}}{\delta t}=\left(\mathbf{A}\right)_{;k}\frac{du^{k}}{dt}
\end{equation}
it can be seen as an instance of the chain rule of differentiation
where the two contracted indices represent the in-between coordinate
differential.

$\bullet$ Because the absolute derivative along a curve is just an
inner product of the covariant derivative with the vector tangent
to the curve, the well known rules of ordinary differentiation of
sums and products also apply to absolute differentiation, as for covariant
differentiation, that is:
\begin{equation}
\frac{\delta}{\delta t}\left(a\mathbf{A}+b\mathbf{B}\right)=a\frac{\delta\mathbf{A}}{\delta t}+b\frac{\delta\mathbf{B}}{\delta t}
\end{equation}
\begin{equation}
\frac{\delta}{\delta t}\left(\mathbf{A}\circ\mathbf{B}\right)=\left(\frac{\delta\mathbf{A}}{\delta t}\circ\mathbf{B}\right)+\left(\mathbf{A}\circ\frac{\delta\mathbf{B}}{\delta t}\right)
\end{equation}
where $a$ and $b$ are constant scalars, $\mathbf{A}$ and $\mathbf{B}$
are differentiable tensors and the symbol $\circ$ denotes an inner
or outer product of tensors.

$\bullet$ The covariant and contravariant metric tensors are in lieu
of constants with respect to absolute differentiation, that is:
\begin{equation}
\frac{\delta g_{ij}}{\delta t}=0\,\,\,\,\,\,\,\,\,\,\,\,\,\,\&\,\,\,\,\,\,\,\,\,\,\,\,\,\,\frac{\delta g^{ij}}{\delta t}=0
\end{equation}
and hence they pass through the absolute derivative operator:
\begin{equation}
\frac{\delta\left(g_{ij}A^{j}\right)}{\delta t}=g_{ij}\frac{\delta A^{j}}{\delta t}\,\,\,\,\,\,\,\,\,\,\,\,\,\,\,\&\,\,\,\,\,\,\,\,\,\,\,\,\,\,\,\,\,\frac{\delta\left(g^{ij}A_{j}\right)}{\delta t}=g^{ij}\frac{\delta A_{j}}{\delta t}
\end{equation}

$\bullet$ For coordinate systems in which all the components of the
metric tensor are constants, the absolute derivative is the same as
the ordinary derivative. This is the case in the rectilinear coordinate
systems.

$\bullet$ The absolute derivative of a tensor along a given curve
is unique, and hence the ordinary derivative of the tensor along that
curve in a rectangular coordinate system is the same as the absolute
derivative of the tensor along that curve in any other system.

\pagebreak{}

\section{Differential Operations}

$\bullet$ In this section we generalize and expand what have been
given in the previous notes \cite{SochiTensorIntro2016} about the
main differential operations which are based on the nabla operator
$\nabla$. The section will investigate these operations in general
curvilinear coordinate systems and in general orthogonal coordinate
systems which are a special case of the general curvilinear systems.
We also investigate the two most important and widely used non-Cartesian
orthogonal coordinate systems, namely the cylindrical and spherical
systems, due to their particular importance.

\subsection{General Curvilinear Coordinate System\label{subGeneralCurvilinearSystem}}

$\bullet$ Here, we investigate the differential operations and operators
in general curvilinear coordinate systems, whether orthogonal or not.

$\bullet$ The previous definitions of the differential operations,
as given in the first set of notes, are essentially valid in general
non-Cartesian coordinate systems if the operations are extended to
include the basis vectors as well as the components.

$\bullet$ The analytical expressions of the differential operations
can be obtained directly if the expression for the nabla operator
$\nabla$ and the spatial derivatives of the basis vectors in the
general curvilinear coordinate system are known.

\subsubsection{Gradient}

$\bullet$ The nabla operator $\nabla$ in general curvilinear coordinate
systems is defined as follow:
\begin{equation}
\nabla=\mathbf{E}^{i}\partial_{i}
\end{equation}
Hence, the gradient of a differentiable scalar function of position,
$f$, is given by:
\begin{equation}
\nabla f=\mathbf{E}^{i}\partial_{i}f=\mathbf{E}^{i}f_{,i}\label{eqGradientGeneral}
\end{equation}
The components of this expression represent the covariant form of
a rank-1 tensor, i.e. $\left[\nabla f\right]_{i}=f_{,i}$, as it should
be since the gradient operation increases the covariant rank of a
tensor by one. Since this expression consists of a contravariant basis
vector and a covariant component, the gradient in general curvilinear
systems is invariant under admissible transformations of coordinates.

$\bullet$ The contravariant form of the gradient of a scalar $f$
can be obtained by using the index raising operator, that is:
\begin{equation}
\left[\nabla f\right]^{i}=\partial^{i}f=g^{ij}\partial_{j}f=g^{ij}f_{,j}=f^{,i}
\end{equation}

$\bullet$ The gradient of a differentiable covariant vector $\mathbf{A}$
can similarly be defined as follow:
\begin{equation}
\begin{aligned}\nabla\mathbf{A} & =\mathbf{E}^{i}\partial_{i}\left(A_{j}\mathbf{E}^{j}\right) & \,\,\,\,\,\,\,\,\,\,\,\,\,\\
 & =\mathbf{E}^{i}\mathbf{E}^{j}\partial_{i}A_{j}+\mathbf{E}^{i}A_{j}\partial_{i}\mathbf{E}^{j} &  & \text{(product rule)}\\
 & =\mathbf{E}^{i}\mathbf{E}^{j}\partial_{i}A_{j}+\mathbf{E}^{i}A_{j}\left(-\Gamma_{ki}^{j}\mathbf{E}^{k}\right) &  & \text{(Eq. \ref{eqPartialBasisContravariant})}\\
 & =\mathbf{E}^{i}\mathbf{E}^{j}\partial_{i}A_{j}-\mathbf{E}^{i}\mathbf{E}^{j}\Gamma_{ji}^{k}A_{k} &  & \text{(relabeling dummy indices \ensuremath{j} \& \ensuremath{k})}\\
 & =\mathbf{E}^{i}\mathbf{E}^{j}\left(\partial_{i}A_{j}-\Gamma_{ji}^{k}A_{k}\right) &  & \text{(taking common factor \ensuremath{\mathbf{E}^{i}\mathbf{E}^{j}})}\\
 & =\mathbf{E}^{i}\mathbf{E}^{j}A_{j;i} &  & \text{(definition of covariant derivative)}
\end{aligned}
\end{equation}
Similarly, for a differentiable contravariant vector $\mathbf{A}$
the gradient is given by:
\begin{equation}
\begin{aligned}\nabla\mathbf{A} & =\mathbf{E}^{i}\partial_{i}\left(A^{j}\mathbf{E}_{j}\right) & \,\,\,\,\,\,\,\,\,\,\,\,\,\,\,\,\\
 & =\mathbf{E}^{i}\mathbf{E}_{j}\partial_{i}A^{j}+\mathbf{E}^{i}A^{j}\partial_{i}\mathbf{E}_{j} &  & \text{(product rule)}\\
 & =\mathbf{E}^{i}\mathbf{E}_{j}\partial_{i}A^{j}+\mathbf{E}^{i}A^{j}\left(\Gamma_{ji}^{k}\mathbf{E}_{k}\right) &  & \text{(Eq. \ref{eqPartialBasisCovariant})}\\
 & =\mathbf{E}^{i}\mathbf{E}_{j}\partial_{i}A^{j}+\mathbf{E}^{i}\mathbf{E}_{j}\Gamma_{ki}^{j}A^{k} &  & \text{(relabeling dummy indices \ensuremath{j} \& \ensuremath{k})}\\
 & =\mathbf{E}^{i}\mathbf{E}_{j}\left(\partial_{i}A^{j}+\Gamma_{ki}^{j}A^{k}\right) &  & \text{(taking common factor \ensuremath{\mathbf{E}^{i}\mathbf{E}_{j}})}\\
 & =\mathbf{E}^{i}\mathbf{E}_{j}A_{\,\,;i}^{j} &  & \text{(definition of covariant derivative)}
\end{aligned}
\end{equation}
The components of the gradients of covariant and contravariant vectors
represent, respectively, the covariant and mixed forms of a rank-2
tensor, as they should be since the gradient operation increases the
covariant rank of a tensor by one.

$\bullet$ The gradient of higher rank tensors is similarly defined.
For example, the gradient of a rank-2 tensor is given by:
\begin{eqnarray}
\nabla\mathbf{A} & = & \mathbf{E}^{i}\mathbf{E}^{j}\mathbf{E}^{k}\left(\partial_{i}A_{jk}-\Gamma_{ji}^{l}A_{lk}-\Gamma_{ki}^{l}A_{jl}\right)=\mathbf{E}^{i}\mathbf{E}^{j}\mathbf{E}^{k}A_{jk;i}\,\,\,\,\,\,\,\,\,\,\,\,\text{(covariant)}\\
\nabla\mathbf{A} & = & \mathbf{E}^{i}\mathbf{E}_{j}\mathbf{E}_{k}\left(\partial_{i}A^{jk}+\Gamma_{li}^{j}A^{lk}+\Gamma_{li}^{k}A^{jl}\right)=\mathbf{E}^{i}\mathbf{E}_{j}\mathbf{E}_{k}A_{\,;i}^{jk}\,\,\,\,\,\,\,\,\,\,\,\,\,\,\,\,\text{(contravariant)}\\
\nabla\mathbf{A} & = & \mathbf{E}^{i}\mathbf{E}^{j}\mathbf{E}_{k}\left(\partial_{i}A_{j}^{\,\,k}-\Gamma_{ji}^{l}A_{l}^{\,\,k}+\Gamma_{li}^{k}A_{j}^{\,\,l}\right)=\mathbf{E}^{i}\mathbf{E}^{j}\mathbf{E}_{k}A_{j;i}^{\,\,k}\,\,\,\,\,\,\,\,\,\,\,\,\,\,\,\,\text{(mixed)}\\
\nabla\mathbf{A} & = & \mathbf{E}^{i}\mathbf{E}_{j}\mathbf{E}^{k}\left(\partial_{i}A_{\,\,k}^{j}+\Gamma_{li}^{j}A_{\,\,k}^{l}-\Gamma_{ki}^{l}A_{\,\,l}^{j}\right)=\mathbf{E}^{i}\mathbf{E}_{j}\mathbf{E}^{k}A_{\,\,k;i}^{j}\,\,\,\,\,\,\,\,\,\,\,\,\,\,\,\text{(mixed)}
\end{eqnarray}

\subsubsection{Divergence}

$\bullet$ Generically, the divergence of a differentiable contravariant
vector $\mathbf{A}$ is defined as follow:
\begin{equation}
\nabla\cdot\mathbf{A}=\mathbf{E}^{i}\partial_{i}\cdot\left(A^{j}\mathbf{E}_{j}\right)=\mathbf{E}^{i}\cdot\partial_{i}\left(A^{j}\mathbf{E}_{j}\right)=\mathbf{E}^{i}\cdot\left(A_{;i}^{j}\mathbf{E}_{j}\right)=\left(\mathbf{E}^{i}\cdot\mathbf{E}_{j}\right)A_{\,;i}^{j}=\delta_{j}^{i}A_{\,;i}^{j}=A_{\,;i}^{i}
\end{equation}
In more details, the divergence of a differentiable contravariant
vector $A^{i}$ is a scalar obtained by contracting the covariant
derivative index with the contravariant index of the vector, and hence:
\begin{equation}
\begin{aligned}\nabla\cdot\mathbf{A} & =A_{\,;i}^{i} & \,\,\,\,\,\,\,\,\,\,\,\,\,\,\,\,\\
 & =\partial_{i}A^{i}+\Gamma_{ji}^{i}A^{j} &  & \text{(definition of covariant derivative)}\\
 & =\partial_{i}A^{i}+A^{j}\frac{1}{\sqrt{g}}\partial_{j}\left(\sqrt{g}\right) &  & \text{(Eq. \ref{eqGammaaia})}\\
 & =\partial_{i}A^{i}+A^{i}\frac{1}{\sqrt{g}}\partial_{i}\left(\sqrt{g}\right) &  & \text{(renaming dummy index \ensuremath{j})}\\
 & =\frac{1}{\sqrt{g}}\partial_{i}\left(\sqrt{g}A^{i}\right) &  & \text{(product rule)}
\end{aligned}
\label{eqDivGeneral}
\end{equation}
where $g$ is the determinant of the covariant metric tensor $g_{ij}$.
The last equality may be called the Voss-Weyl formula.

$\bullet$ The divergence can also be obtained by raising the first
index of the covariant derivative of a covariant vector using a contracting
contravariant metric tensor:
\begin{equation}
\begin{aligned}g^{ji}A_{j;i} & =\left(g^{ji}A_{j}\right)_{;i} & \,\,\,\,\,\,\,\,\,\,\,\,\,\,\, & \text{(Eq. \ref{eqCovDerivMetric})}\\
 & =\left(A^{i}\right)_{;i} &  & \text{(raising index)}\\
 & =A_{\,\,;i}^{i}\\
 & =\nabla\cdot\mathbf{A} &  & \text{(Eq. \ref{eqDivGeneral})}
\end{aligned}
\end{equation}
as before.

$\bullet$ Based on the previous point, the divergence of a covariant
vector $A_{j}$ is obtained by using the raising operator, that is
\begin{equation}
A_{\,;i}^{i}=g^{ij}A_{j;i}
\end{equation}

$\bullet$ For a rank-2 contravariant tensor $\mathbf{A}$, the divergence
is generically defined by:
\begin{equation}
\nabla\cdot\mathbf{A}=\mathbf{E}^{i}\partial_{i}\cdot\left(A^{jk}\mathbf{E}_{j}\mathbf{E}_{k}\right)=\left(\mathbf{E}^{i}\cdot\mathbf{E}_{j}\right)\mathbf{E}_{k}A_{\,\,;i}^{jk}=\delta_{j}^{i}\mathbf{E}_{k}A_{\,\,;i}^{jk}=\mathbf{E}_{k}A_{\,\,;i}^{ik}
\end{equation}
The components of this expression represent a contravariant vector,
as it should be since the divergence operation reduces the contravariant
rank of a tensor by one.

$\bullet$ More generally, considering the tensor components, the
divergence of a differentiable rank-2 contravariant tensor $A^{ij}$
is a contravariant vector obtained by contracting the covariant derivative
index with one of the contravariant indices, e.g.
\begin{equation}
\left[\nabla\cdot\mathbf{A}\right]^{j}=A_{\,\,;i}^{ij}
\end{equation}
And for a rank-2 mixed tensor $A_{j}^{i}$ we have:
\begin{equation}
\left[\nabla\cdot\mathbf{A}\right]_{j}=A_{j;i}^{i}
\end{equation}

$\bullet$ Similarly, for a general tensor of type ($m,n$): $\mathbf{A}=A_{j_{1}j_{2}\cdots\cdots\cdots j_{n}}^{i_{1}i_{2}\cdots i_{k}\cdots i_{m}}$,
the divergence with respect to its $k^{th}$ contravariant index is
defined by:
\begin{equation}
\left[\nabla\cdot\mathbf{A}\right]_{j_{1}j_{2}\cdots\cdots\cdots j_{n}}^{i_{1}i_{2}\cdots\,\,\cdots i_{m}}=\left(A_{j_{1}j_{2}\cdots\cdots\cdots j_{n}}^{i_{1}i_{2}\cdots s\cdots i_{m}}\right)_{;s}
\end{equation}
with the absence of the contracted contravariant index $i_{k}$ on
the left hand side.

\subsubsection{Curl}

$\bullet$ The curl of a differentiable vector is the cross product
of the nabla operator $\nabla$ with the vector, e.g. the curl of
a vector $\mathbf{A}$ represented by covariant components is given
by:
\begin{equation}
\begin{aligned}\mathrm{curl\,}\mathbf{A} & =\nabla\times\mathbf{A} & \,\,\,\\
 & =\mathbf{E}^{i}\partial_{i}\times A_{j}\mathbf{E}^{j}\\
 & =\mathbf{E}^{i}\times\partial_{i}\left(A_{j}\mathbf{E}^{j}\right)\\
 & =\mathbf{E}^{i}\times\left(A_{j;i}\mathbf{E}^{j}\right) &  & \text{(Eq. \ref{eqCoDifCo})}\\
 & =A_{j;i}\left(\mathbf{E}^{i}\times\mathbf{E}^{j}\right)\\
 & =A_{j;i}\,\underline{\epsilon}^{ijk}\mathbf{E}_{k} &  & \text{(Eq. \ref{eqEiEjCon})}\\
 & =\underline{\epsilon}^{ijk}A_{j;i}\mathbf{E}_{k}\\
 & =\frac{\epsilon^{ijk}}{\sqrt{g}}\left(\partial_{i}A_{j}-\Gamma_{ji}^{l}A_{l}\right)\mathbf{E}_{k} &  & \text{(Eq. \ref{eqepsilonunderline} \& definition of covariant derivative)}
\end{aligned}
\end{equation}
and hence the $k^{th}$ contravariant component of $\mathrm{curl\,}\mathbf{A}$
is:
\begin{equation}
\left[\nabla\times\mathbf{A}\right]^{k}=\frac{\epsilon^{ijk}}{\sqrt{g}}\left(\partial_{i}A_{j}-\Gamma_{ji}^{l}A_{l}\right)\label{eqCurlGeneralComponent}
\end{equation}
On expanding the last equation for the three components of a 3D space,
considering that the terms of the Christoffel symbols cancel out due
to their symmetry in the two lower indices,\footnote{That is: $A_{i;j}-A_{j;i}=\partial_{j}A_{i}-A_{k}\Gamma_{ij}^{k}-\partial_{i}A_{j}+A_{k}\Gamma_{ji}^{k}=\partial_{j}A_{i}-A_{k}\Gamma_{ij}^{k}-\partial_{i}A_{j}+A_{k}\Gamma_{ij}^{k}=\partial_{j}A_{i}-\partial_{i}A_{j}=A_{i,j}-A_{j,i}$.}
we obtain:
\begin{eqnarray}
\left[\nabla\times\mathbf{A}\right]^{1} & = & \frac{1}{\sqrt{g}}\left(\partial_{2}A_{3}-\partial_{3}A_{2}\right)\\
\left[\nabla\times\mathbf{A}\right]^{2} & = & \frac{1}{\sqrt{g}}\left(\partial_{3}A_{1}-\partial_{1}A_{3}\right)\\
\left[\nabla\times\mathbf{A}\right]^{3} & = & \frac{1}{\sqrt{g}}\left(\partial_{1}A_{2}-\partial_{2}A_{1}\right)
\end{eqnarray}
Hence, Eq. \ref{eqCurlGeneralComponent} will reduce to:
\begin{equation}
\left[\nabla\times\mathbf{A}\right]^{k}=\frac{\epsilon^{ijk}}{\sqrt{g}}\partial_{i}A_{j}
\end{equation}

\subsubsection{Laplacian\label{subLaplacianGneralCurvilinear}}

$\bullet$ Generically, the Laplacian of a differentiable scalar function
of position, $f$, is defined as follow:
\begin{equation}
\nabla^{2}f=\mathrm{div}\left(\mathrm{grad\,}f\right)=\nabla\cdot\left(\nabla f\right)
\end{equation}
Hence the simplest approach for obtaining the Laplacian in general
coordinates is to insert the expression for the gradient, $\nabla f$,
into the expression for the divergence. However, because in general
curvilinear coordinates the divergence is defined only for contravariant
tensors whereas the gradient of a scalar is a covariant tensor, the
index of the gradient should be raised first before applying the divergence
operation, that is:
\begin{equation}
\left[\nabla f\right]^{i}=\partial^{i}f=g^{ij}\partial_{j}f
\end{equation}
Now, according to Eq. \ref{eqDivGeneral} the divergence is given
by:
\begin{equation}
\nabla\cdot\mathbf{A}=A_{\,\,;i}^{i}=\frac{1}{\sqrt{g}}\partial_{i}\left(\sqrt{g}A^{i}\right)
\end{equation}
On defining $\mathbf{A}\equiv\nabla f=\mathbf{E}_{i}\partial^{i}f$
and replacing $A^{i}$ in the last equation with $\partial^{i}f$
we obtain:
\begin{equation}
\nabla^{2}f=\nabla\cdot\left(\mathbf{E}_{i}\partial^{i}f\right)=\frac{1}{\sqrt{g}}\partial_{i}\left(\sqrt{g}g^{ij}\partial_{j}f\right)\label{eqLaplacianGeneral}
\end{equation}
which is the expression for the Laplacian of a scalar function $f$
in general curvilinear coordinates.

$\bullet$ Another approach for developing the Laplacian expression
in general coordinates is to apply the first principles by using the
definitions and basic properties of the operations involved, that
is:
\begin{equation}
\begin{aligned}\nabla^{2}f & =\nabla\cdot\left(\nabla f\right)\\
 & =\mathbf{E}^{i}\partial_{i}\cdot\left(\mathbf{E}^{j}\partial_{j}f\right)\\
 & =\mathbf{E}^{i}\cdot\partial_{i}\left(\mathbf{E}^{j}\partial_{j}f\right)\\
 & =\mathbf{E}^{i}\cdot\partial_{i}\left(\mathbf{E}^{j}f_{,j}\right)\\
 & =\mathbf{E}^{i}\cdot\left(\mathbf{E}^{j}f_{,j;i}\right) &  & \text{(Eq. \ref{eqCoDifCo})}\\
 & =\left(\mathbf{E}^{i}\cdot\mathbf{E}^{j}\right)f_{,j;i}\\
 & =g^{ij}f_{,j;i} &  & \text{(Eq. \ref{eqEg})}\\
 & =\left(g^{ij}f_{,j}\right)_{;i} &  & \text{(Eq. \ref{eqCovDerivMetric})}\\
 & =\left(g^{ij}\partial_{j}f\right)_{;i}\\
 & =\partial_{i}\left(g^{ij}\partial_{j}f\right)+\left(g^{kj}\partial_{j}f\right)\Gamma_{ki}^{i} &  & \text{(Eq. \ref{eqCovariantDerivative})}\\
 & =\partial_{i}\left(g^{ij}\partial_{j}f\right)+\left(g^{ij}\partial_{j}f\right)\Gamma_{ik}^{k} &  & \text{(renaming dummy indices \ensuremath{i} \& \ensuremath{k})}\\
 & =\partial_{i}\left(g^{ij}\partial_{j}f\right)+\left(g^{ij}\partial_{j}f\right)\frac{1}{\sqrt{g}}\left(\partial_{i}\sqrt{g}\right) &  & \text{(Eq. \ref{eqGammaaia})}\\
 & =\frac{1}{\sqrt{g}}\left[\sqrt{g}\partial_{i}\left(g^{ij}\partial_{j}f\right)+\left(g^{ij}\partial_{j}f\right)\partial_{i}\sqrt{g}\right] &  & \text{(taking \ensuremath{\frac{1}{\sqrt{g}}} factor)}\\
 & =\frac{1}{\sqrt{g}}\partial_{i}\left(\sqrt{g}g^{ij}\partial_{j}f\right) &  & \text{(product rule)}
\end{aligned}
\end{equation}
as before.

$\bullet$ The Laplacian of a scalar $f$ may also be shorthand notated
with:
\begin{equation}
\nabla^{2}f=g^{ij}f_{,ij}
\end{equation}

$\bullet$ The Laplacian of non-scalar tensors can be similarly defined.
For example, the Laplacian of a vector $\mathbf{B}$ is a vector $\mathbf{A}$
(i.e. $\mathbf{A}=\nabla^{2}\mathbf{B}$) which may be defined in
general coordinates as:
\begin{equation}
A^{i}=g^{jk}B_{\,\,\,;jk}^{i}\,\,\,\,\,\,\,\,\,\,\,\,\text{and}\,\,\,\,\,\,\,\,\,\,\,A_{i}=g^{jk}B_{i;jk}
\end{equation}

$\bullet$ The Laplacian of a tensor is a tensor of the same rank
and variance type.

\subsection{General Orthogonal Coordinate System\label{subOrthogonalCurvilinearSystem}}

$\bullet$ In this section we state the main differential operations
in general orthogonal coordinate systems. These operations are special
cases of the operations in general curvilinear systems which were
derived in the previous section. However, due to the wide spread use
of orthogonal systems, it is worth to state the most important of
these operations although they can be easily obtained from the formulae
of general curvilinear systems.

$\bullet$ General orthogonal coordinate systems are identified in
the following notes by the coordinates ($u^{1},\ldots,u^{n}$) with
unit basis vectors ($\mathbf{u}_{1},\ldots,\mathbf{u}_{n}$) and scale
factors ($h_{1},\ldots,h_{n}$) where:\footnote{In orthogonal coordinates, the covariant and contravariant normalized
basis vectors are identical, as established previously in $\S$ \ref{subCovariantContravariantPhysical},
and hence $\mathbf{u}^{i}=\mathbf{u}_{i}$ and $\mathbf{e}^{j}=\mathbf{e}_{j}$.}
\begin{equation}
\mathbf{u}_{i}=\sum_{j}\frac{1}{h_{i}}\frac{\partial x^{j}}{\partial u^{i}}\mathbf{e}_{j}=\sum_{j}h_{i}\frac{\partial u^{i}}{\partial x^{j}}\mathbf{e}_{j}\,\,\,\,\,\,\,\,\,\,\,\,\,\,\,\text{(no sum on \ensuremath{i})}
\end{equation}
\begin{equation}
h_{i}=\left|\mathbf{E}_{i}\right|=\left|\frac{\partial\mathbf{r}}{\partial u_{i}}\right|=\left[\sum_{j}\left(\frac{\partial x^{j}}{\partial u^{i}}\right)^{2}\right]^{1/2}=\left[\sum_{j}\left(\frac{\partial u^{i}}{\partial x^{j}}\right)^{2}\right]^{-1/2}
\end{equation}
In the last equations, $x^{j}$ and $\mathbf{e}_{j}$ are respectively
the coordinates and unit basis vectors in the Cartesian rectangular
system, and $\mathbf{r}$ is the position vector in that system.

\subsubsection{Gradient}

$\bullet$ The nabla operator in general orthogonal coordinates is
given by:
\begin{equation}
\nabla=\sum_{i}\frac{\mathbf{u}_{i}}{h_{i}}\frac{\partial}{\partial u^{i}}
\end{equation}
Hence, the gradient of a differentiable scalar $f$ in orthogonal
coordinates is given by:
\begin{equation}
\nabla f=\sum_{i}\frac{\mathbf{u}_{i}}{h_{i}}\frac{\partial f}{\partial u^{i}}
\end{equation}

\subsubsection{Divergence}

$\bullet$ The divergence in orthogonal coordinates can be obtained
from Eq. \ref{eqDivGeneral}. Since for orthogonal coordinate systems
the metric tensor is diagonal with $\sqrt{g}=h_{1}h_{2}h_{3}$ in
a 3D space and $h_{i}A^{i}=\hat{A}^{i}$ (component-wise with no summation),
the last line of Eq. \ref{eqDivGeneral} becomes:
\begin{equation}
\nabla\cdot\mathbf{A}=\frac{1}{\sqrt{g}}\frac{\partial}{\partial u^{i}}\left(\sqrt{g}A^{i}\right)=\frac{1}{h_{1}h_{2}h_{3}}\sum_{i=1}^{3}\frac{\partial}{\partial u^{i}}\left(\frac{h_{1}h_{2}h_{3}}{h_{i}}\hat{A}^{i}\right)
\end{equation}
where $\mathbf{A}$ is a contravariant differentiable vector and $\hat{A}^{i}$
represents the physical components (refer to Eq. \ref{eqPhysCompCova}).\footnote{In orthogonal coordinate systems the physical components are the same
for covariant and contravariant forms, as established before in $\S$
\ref{subCovariantContravariantPhysical}, and hence $\hat{A}_{i}=\hat{A}^{i}$.}~This equation is the divergence of a vector in general orthogonal
coordinates as defined in vector calculus.

\subsubsection{Curl}

$\bullet$ The curl of a differentiable vector $\mathbf{A}$ in orthogonal
coordinate systems in 3D spaces is given by:
\begin{equation}
\nabla\times\mathbf{A}=\frac{1}{h_{1}h_{2}h_{3}}\begin{vmatrix}\begin{array}{ccc}
h_{1}\mathbf{u}_{1} & h_{2}\mathbf{u}_{2} & h_{3}\mathbf{u}_{3}\\
\frac{\partial}{\partial u^{1}} & \frac{\partial}{\partial u^{2}} & \frac{\partial}{\partial u^{3}}\\
h_{1}\hat{A}_{1} & h_{2}\hat{A}_{2} & h_{3}\hat{A}_{3}
\end{array}\end{vmatrix}
\end{equation}
where the hat indicates a physical component. The last equation may
also be given in a more compact form as:
\begin{equation}
\left[\nabla\times\mathbf{A}\right]_{i}=\sum_{k=1}^{3}\frac{\epsilon_{ijk}h_{i}}{h_{1}h_{2}h_{3}}\frac{\partial(h_{k}\hat{A}_{k})}{\partial u^{j}}
\end{equation}

\subsubsection{Laplacian}

$\bullet$ For general orthogonal coordinate systems in 3D spaces
we have:
\begin{equation}
\sqrt{g}=h_{1}h_{2}h_{3}\,\,\,\,\,\,\,\,\,\,\,\,\,\&\,\,\,\,\,\,\,\,\,\,\,\,\,g^{ij}=\frac{\delta^{ij}}{h_{i}h_{j}}\,\,\,\,\,\,\,\text{(no sum on \ensuremath{i} or \ensuremath{j})}
\end{equation}
and hence Eq. \ref{eqLaplacianGeneral} becomes:
\begin{equation}
\nabla^{2}f=\frac{1}{h_{1}h_{2}h_{3}}\sum_{i=1}^{3}\frac{\partial}{\partial u^{i}}\left(\frac{h_{1}h_{2}h_{3}}{h_{i}^{2}}\frac{\partial f}{\partial u^{i}}\right)
\end{equation}
which is the Laplacian of a scalar function of position, $f$, in
orthogonal coordinates as defined in vector calculus.

\subsection{Cylindrical Coordinate System\label{subCylindricalSystem}}

$\bullet$ For cylindrical systems identified by the coordinates ($\rho,\phi,z$),
the orthonormal basis vectors are $\mathbf{e}_{\rho},\,\mathbf{e}_{\phi}$
and $\mathbf{e}_{z}$.\footnote{Hence the given components (i.e. $A_{\rho},A_{\phi}$ and $A_{z}$)
are physical (see $\S$ \ref{subCovariantContravariantPhysical}).
Despite that, we do not use hats since the components are suffixed
with coordinate symbols (refer to $\S$ \ref{subCovariantContravariantPhysical}).}~We use for brevity $\mathbf{e}_{\rho\phi}$ as a shorthand notation
for the unit dyad $\mathbf{e}_{\rho}\mathbf{e}_{\phi}$ and similar
notations for the other dyads.

\subsubsection{Gradient}

$\bullet$ The gradient of a differentiable scalar $f$ is:
\begin{equation}
\nabla f=\mathbf{e}_{\rho}\partial_{\rho}f+\mathbf{e}_{\phi}\frac{1}{\rho}\partial_{\phi}f+\mathbf{e}_{z}\partial_{z}f
\end{equation}

$\bullet$ The gradient of a differentiable vector $\mathbf{A}$ is:
\begin{eqnarray}
\nabla\mathbf{A} & = & \mathbf{e}_{\rho\rho}A_{\rho,\rho}+\mathbf{e}_{\rho\phi}A_{\phi,\rho}+\mathbf{e}_{\rho z}A_{z,\rho}+\label{eqGradVectorCylindrical}\\
 &  & \mathbf{e}_{\phi\rho}\left(\frac{1}{\rho}A_{\rho,\phi}-\frac{A_{\phi}}{\rho}\right)+\mathbf{e}_{\phi\phi}\left(\frac{1}{\rho}A_{\phi,\phi}+\frac{A_{\rho}}{\rho}\right)+\mathbf{e}_{\phi z}\frac{1}{\rho}A_{z,\phi}+\nonumber \\
 &  & \mathbf{e}_{z\rho}A_{\rho,z}+\mathbf{e}_{z\phi}A_{\phi,z}+\mathbf{e}_{zz}A_{z,z}\nonumber
\end{eqnarray}

\subsubsection{Divergence}

$\bullet$ The divergence of a differentiable vector $\mathbf{A}$
is:
\begin{equation}
\nabla\cdot\mathbf{A}=\frac{1}{\rho}\left[\partial_{\rho}\left(\rho A_{\rho}\right)+\partial_{\phi}A_{\phi}+\rho\partial_{z}A_{z}\right]
\end{equation}

$\bullet$ The divergence of a differentiable rank-2 tensor $\mathbf{A}$
is a vector given by:\footnote{It should be understood that $\nabla\cdot\mathbf{A}$ is lower than
the original tensor $\mathbf{A}$ by just one contravariant index
and hence, unlike the common use of this notation, it is not scalar
in general.}
\begin{eqnarray}
\nabla\cdot\mathbf{A} & = & \mathbf{e}_{\rho}\left(A_{\rho\rho,\rho}+\frac{A_{\rho\rho}-A_{\phi\phi}}{\rho}+\frac{1}{\rho}A_{\phi\rho,\phi}+A_{z\rho,z}\right)+\label{eqDivTensorCylindrical}\\
 &  & \mathbf{e}_{\phi}\left(A_{\rho\phi,\rho}+\frac{2A_{\rho\phi}}{\rho}+\frac{1}{\rho}A_{\phi\phi,\phi}+A_{z\phi,z}+\frac{A_{\phi\rho}-A_{\rho\phi}}{\rho}\right)+\nonumber \\
 &  & \mathbf{e}_{z}\left(A_{\rho z,\rho}+\frac{A_{\rho z}}{\rho}+\frac{1}{\rho}A_{\phi z,\phi}+A_{zz,z}\right)\nonumber
\end{eqnarray}

\subsubsection{Curl}

$\bullet$ The curl of a differentiable vector $\mathbf{A}$ is:
\begin{equation}
\nabla\times\mathbf{A}=\frac{1}{\rho}\begin{vmatrix}\begin{array}{ccc}
\mathbf{e}_{\rho} & \rho\mathbf{e}_{\phi} & \mathbf{e}_{z}\\
\partial_{\rho} & \partial_{\phi} & \partial_{z}\\
A_{\rho} & \rho A_{\phi} & A_{z}
\end{array}\end{vmatrix}
\end{equation}

\subsubsection{Laplacian}

$\bullet$ The Laplacian of a differentiable scalar $f$ is:
\begin{equation}
\nabla^{2}f=\partial_{\rho\rho}f+\frac{1}{\rho}\partial_{\rho}f+\frac{1}{\rho^{2}}\partial_{\phi\phi}f+\partial_{zz}f
\end{equation}

$\bullet$ The Laplacian of a differentiable vector $\mathbf{A}$
is:
\begin{eqnarray}
\nabla^{2}\mathbf{A} & = & \mathbf{e}_{\rho}\left[\partial_{\rho}\left(\frac{1}{\rho}\partial_{\rho}\left(\rho A_{\rho}\right)\right)+\frac{1}{\rho^{2}}\partial_{\phi\phi}A_{\rho}+\partial_{zz}A_{\rho}-\frac{2}{\rho^{2}}\partial_{\phi}A_{\phi}\right]+\label{eqLaplacianVectorCylindrical}\\
 &  & \mathbf{e}_{\phi}\left[\partial_{\rho}\left(\frac{1}{\rho}\partial_{\rho}\left(\rho A_{\phi}\right)\right)+\frac{1}{\rho^{2}}\partial_{\phi\phi}A_{\phi}+\partial_{zz}A_{\phi}+\frac{2}{\rho^{2}}\partial_{\phi}A_{\rho}\right]+\nonumber \\
 &  & \mathbf{e}_{z}\left[\frac{1}{\rho}\partial_{\rho}\left(\rho\partial_{\rho}A_{z}\right)+\frac{1}{\rho^{2}}\partial_{\phi\phi}A_{z}+\partial_{zz}A_{z}\right]\nonumber
\end{eqnarray}

\subsection{Spherical Coordinate System\label{subSphericalSystem}}

$\bullet$ For spherical coordinate systems identified by the coordinates
($r,\theta,\phi$), the orthonormal basis vectors are $\mathbf{e}_{r},\,\mathbf{e}_{\theta}$
and $\mathbf{e}_{\phi}$.\footnote{Again, the components are physical and we do not use hats.}
We use for brevity $\mathbf{e}_{r\theta}$ as a shorthand notation
for the dyad $\mathbf{e}_{r}\mathbf{e}_{\theta}$ and similar notations
for the other unit dyads.

\subsubsection{Gradient}

$\bullet$ The gradient of a differentiable scalar $f$ is:
\begin{equation}
\nabla f=\mathbf{e}_{r}\partial_{r}f+\mathbf{e}_{\theta}\frac{1}{r}\partial_{\theta}f+\mathbf{e}_{\phi}\frac{1}{r\sin\theta}\partial_{\phi}f
\end{equation}

$\bullet$ The gradient of a differentiable vector $\mathbf{A}$ is:
\begin{eqnarray}
\nabla\mathbf{A} & = & \mathbf{e}_{rr}A_{r,r}+\mathbf{e}_{r\theta}A_{\theta,r}+\mathbf{e}_{r\phi}A_{\phi,r}+\label{eqGradVectorSpherical}\\
 &  & \mathbf{e}_{\theta r}\left(\frac{A_{r,\theta}}{r}-\frac{A_{\theta}}{r}\right)+\mathbf{e}_{\theta\theta}\left(\frac{A_{\theta,\theta}}{r}+\frac{A_{r}}{r}\right)+\mathbf{e}_{\theta\phi}\frac{A_{\phi,\theta}}{r}+\nonumber \\
 &  & \mathbf{e}_{\phi r}\left(\frac{A_{r,\phi}}{r\sin\theta}-\frac{A_{\phi}}{r}\right)+\mathbf{e}_{\phi\theta}\left(\frac{A_{\theta,\phi}}{r\sin\theta}-\frac{A_{\phi}\cot\theta}{r}\right)+\mathbf{e}_{\phi\phi}\left(\frac{A_{\phi,\phi}}{r\sin\theta}+\frac{A_{r}}{r}+\frac{A_{\theta}\cot\theta}{r}\right)\nonumber
\end{eqnarray}

\subsubsection{Divergence}

$\bullet$ The divergence of a differentiable vector $\mathbf{A}$
is:
\begin{equation}
\nabla\cdot\mathbf{A}=\frac{1}{r^{2}\sin\theta}\left[\sin\theta\frac{\partial\left(r^{2}A_{r}\right)}{\partial r}+r\frac{\partial\left(\sin\theta A_{\theta}\right)}{\partial\theta}+r\frac{\partial A_{\phi}}{\partial\phi}\right]
\end{equation}

$\bullet$ The divergence of a differentiable rank-2 tensor $\mathbf{A}$
is a vector given by:
\begin{eqnarray}
\nabla\cdot\mathbf{A} & = & \mathbf{e}_{r}\left(\frac{\partial_{r}\left(r^{2}A_{rr}\right)}{r^{2}}+\frac{\partial_{\theta}\left(A_{\theta r}\sin\theta\right)}{r\sin\theta}+\frac{\partial_{\phi}A_{\phi r}}{r\sin\theta}-\frac{A_{\theta\theta}+A_{\phi\phi}}{r}\right)+\label{eqDivTensorSpherical}\\
 &  & \mathbf{e}_{\theta}\left(\frac{\partial_{r}\left(r^{3}A_{r\theta}\right)}{r^{3}}+\frac{\partial_{\theta}\left(A_{\theta\theta}\sin\theta\right)}{r\sin\theta}+\frac{\partial_{\phi}A_{\phi\theta}}{r\sin\theta}+\frac{A_{\theta r}-A_{r\theta}-A_{\phi\phi}\cot\theta}{r}\right)+\nonumber \\
 &  & \mathbf{e}_{\phi}\left(\frac{\partial_{r}\left(r^{3}A_{r\phi}\right)}{r^{3}}+\frac{\partial_{\theta}\left(A_{\theta\phi}\sin\theta\right)}{r\sin\theta}+\frac{\partial_{\phi}A_{\phi\phi}}{r\sin\theta}+\frac{A_{\phi r}-A_{r\phi}+A_{\phi\theta}\cot\theta}{r}\right)\nonumber
\end{eqnarray}

\subsubsection{Curl}

$\bullet$ The curl of a differentiable vector $\mathbf{A}$ is:
\begin{equation}
\nabla\times\mathbf{A}=\frac{1}{r^{2}\sin\theta}\begin{vmatrix}\begin{array}{ccc}
\mathbf{e}_{r} & r\mathbf{e}_{\theta} & r\sin\theta\mathbf{e}_{\phi}\\
\partial_{r} & \partial_{\theta} & \partial_{\phi}\\
A_{r} & rA_{\theta} & r\sin\theta A_{\phi}
\end{array}\end{vmatrix}
\end{equation}

\subsubsection{Laplacian}

$\bullet$ The Laplacian of a differentiable scalar $f$ is:
\begin{equation}
\nabla^{2}f=\partial_{rr}f+\frac{2}{r}\partial_{r}f+\frac{1}{r^{2}}\partial_{\theta\theta}f+\frac{\cos\theta}{r^{2}\sin\theta}\partial_{\theta}f+\frac{1}{r^{2}\sin^{2}\theta}\partial_{\phi\phi}f
\end{equation}

$\bullet$ The Laplacian of a differentiable vector $\mathbf{A}$
is:
\begin{eqnarray}
\hspace{-1cm}\nabla^{2}\mathbf{A} & = & \mathbf{e}_{r}\left[\partial_{r}\left(\frac{\partial_{r}\left(r^{2}A_{r}\right)}{r^{2}}\right)+\frac{\partial_{\theta}\left(\sin\theta\partial_{\theta}A_{r}\right)}{r^{2}\sin\theta}+\frac{\partial_{\phi\phi}A_{r}}{r^{2}\sin^{2}\theta}-\frac{2\partial_{\theta}\left(A_{\theta}\sin\theta\right)}{r^{2}\sin\theta}-\frac{2\partial_{\phi}A_{\phi}}{r^{2}\sin\theta}\right]+\label{eqLaplacianVectorSpherical}\\
 &  & \mathbf{e}_{\theta}\left[\frac{\partial_{r}\left(r^{2}\partial_{r}A_{\theta}\right)}{r^{2}}+\frac{1}{r^{2}}\partial_{\theta}\left(\frac{\partial_{\theta}\left(A_{\theta}\sin\theta\right)}{\sin\theta}\right)+\frac{\partial_{\phi\phi}A_{\theta}}{r^{2}\sin^{2}\theta}+\frac{2\partial_{\theta}A_{r}}{r^{2}}-\frac{2\cot\theta}{r^{2}\sin\theta}\partial_{\phi}A_{\phi}\right]+\nonumber \\
 &  & \mathbf{e}_{\phi}\left[\frac{\partial_{r}\left(r^{2}\partial_{r}A_{\phi}\right)}{r^{2}}+\frac{1}{r^{2}}\partial_{\theta}\left(\frac{\partial_{\theta}\left(A_{\phi}\sin\theta\right)}{\sin\theta}\right)+\frac{\partial_{\phi\phi}A_{\phi}}{r^{2}\sin^{2}\theta}+\frac{2\partial_{\phi}A_{r}}{r^{2}\sin\theta}+\frac{2\cot\theta}{r^{2}\sin\theta}\partial_{\phi}A_{\theta}\right]\nonumber
\end{eqnarray}

\pagebreak{}

\section{Tensors in Applications}

$\bullet$ In this section we conduct a preliminary investigation
of some commonly-used tensors in physical and mathematical applications
of tensor calculus in anticipation of the forthcoming set of notes.
Most of these tensors come from differential geometry, fluid, continuum
and relativistic mechanics, since these disciplines are intimately
linked to tensor calculus as large parts of the subject were developed,
and are still developing, within these disciplines. These tensors
also form the building blocks of several physical and mathematical
theories. We would like to insist that these are just a few partially
representative examples for the use of tensors in scientific and mathematical
applications to have more familiarity with tensor language and techniques
and hence they are not meant to provide a comprehensive view in any
way. We should also indicate that some tensors are defined differently
in different disciplines and hence the given definitions and properties
may not be thorough or general.

\subsection{Riemann Tensor\label{subRiemannTensor}}

$\bullet$ This rank-4 tensor, which is also called Riemann curvature
tensor and Riemann-Christoffel tensor, is a property of the space.
It characterizes important properties of spaces and surfaces and hence
it plays an important role in geometry in general and in non-Euclidean
geometries in particular.

$\bullet$ The covariant differentiation operators in mixed derivatives
are not commutative and hence for a covariant vector $\mathbf{A}$
we have:
\begin{equation}
A_{j;kl}-A_{j;lk}=R_{jkl}^{i}A_{i}\label{eqMixedCovariantDerivative}
\end{equation}
where $R_{jkl}^{i}$ is the Riemann tensor of the second kind which
is given by:
\begin{equation}
R_{\,jkl}^{i}=\partial_{k}\Gamma_{jl}^{i}-\partial_{l}\Gamma_{jk}^{i}+\Gamma_{jl}^{r}\Gamma_{rk}^{i}-\Gamma_{jk}^{r}\Gamma_{rl}^{i}\label{eqRijklMix}
\end{equation}

$\bullet$ The last equation can be put into the following mnemonic
determinantal form:
\begin{equation}
R_{jkl}^{i}=\begin{vmatrix}\begin{array}{cc}
\partial_{k} & \partial_{l}\\
\Gamma_{jk}^{i} & \Gamma_{jl}^{i}
\end{array}\end{vmatrix}+\begin{vmatrix}\begin{array}{cc}
\Gamma_{jl}^{r} & \Gamma_{jk}^{r}\\
\Gamma_{rl}^{i} & \Gamma_{rk}^{i}
\end{array}\end{vmatrix}
\end{equation}

$\bullet$ The Riemann-Christoffel tensor of the second kind is also
called the mixed Riemann-Christoffel tensor.

$\bullet$ From Eq. \ref{eqMixedCovariantDerivative}, it is obvious
that the mixed second order covariant derivatives are equal \textit{iff}
the Riemann tensor of the second kind vanishes identically.

$\bullet$ On lowering the contravariant index of the Riemann tensor
of the second kind, the Riemann tensor of the first kind is obtained:
\begin{equation}
R_{ijkl}=g_{ia}R_{\,\,jkl}^{a}
\end{equation}

$\bullet$ Alternatively, the Riemann tensor of the first kind can
be established independently as:
\begin{eqnarray}
R_{ijkl} & = & \partial_{k}\left[jl,i\right]-\partial_{l}\left[jk,i\right]+\left[il,r\right]\Gamma_{jk}^{r}-\left[ik,r\right]\Gamma_{jl}^{r}\nonumber \\
 & = & \frac{1}{2}\left(\partial_{j}\partial_{k}g_{il}+\partial_{i}\partial_{l}g_{jk}-\partial_{j}\partial_{l}g_{ik}-\partial_{i}\partial_{k}g_{jl}\right)+\left[il,r\right]\Gamma_{jk}^{r}-\left[ik,r\right]\Gamma_{jl}^{r}\\
 & = & \frac{1}{2}\left(\partial_{j}\partial_{k}g_{il}+\partial_{i}\partial_{l}g_{jk}-\partial_{j}\partial_{l}g_{ik}-\partial_{i}\partial_{k}g_{jl}\right)+g^{rs}\left(\left[il,r\right]\left[jk,s\right]-\left[ik,r\right]\left[jl,s\right]\right)\nonumber
\end{eqnarray}

$\bullet$ The first line of the last equation can be cast in the
following mnemonic determinantal form:
\begin{equation}
R_{ijkl}=\begin{vmatrix}\begin{array}{cc}
\partial_{k} & \partial_{l}\\
\left[jk,i\right] & \left[jl,i\right]
\end{array}\end{vmatrix}+\begin{vmatrix}\begin{array}{cc}
\Gamma_{jk}^{r} & \Gamma_{jl}^{r}\\
\left[ik,r\right] & \left[il,r\right]
\end{array}\end{vmatrix}
\end{equation}

$\bullet$ Similarly, the Riemann-Christoffel tensor of the second
kind can be obtained by raising the first covariant index of the Riemann-Christoffel
tensor of the first kind:
\begin{equation}
R_{\,\,\,jkl}^{i}=g^{ia}R_{ajkl}
\end{equation}

$\bullet$ The Riemann-Christoffel tensor of the first kind is also
called the covariant (or totally covariant) Riemann-Christoffel tensor.

$\bullet$ The covariant differentiation operators become commutative
when the metric makes the Riemann tensor of either kind vanish.

$\bullet$ For the mixed second order covariant derivatives of a contravariant
vector $\mathbf{A}$ we have:
\begin{equation}
A_{\,;kl}^{j}-A_{\,;lk}^{j}=R_{\,\,ilk}^{j}A^{i}
\end{equation}
which is similar to Eq. \ref{eqMixedCovariantDerivative} for a covariant
vector $\mathbf{A}$.

$\bullet$ The Riemann-Christoffel tensor vanishes identically \textit{iff}
the space is globally flat. Hence, the Riemann tensor is zero in Euclidean
spaces, and consequently the mixed second order covariant derivatives,
which become ordinary derivatives, are equal when the $C^{2}$ continuity
condition is satisfied.

$\bullet$ The Riemann curvature tensor depends only on the metric
which, in general curvilinear coordinates, is a function of position
and hence the Riemann tensor follows this dependency on position.
Yes, for affine coordinates the metric tensor is constant and hence
the Riemann tensor vanishes identically.

$\bullet$ The totally covariant Riemann tensor satisfies the following
symmetric and skew-symmetric relations in its four indices:
\begin{eqnarray}
R_{ijkl} & = & \,\,\,\,\,R_{klij}\,\,\,\,\,\,\,\,\,\,\,\,\text{(block symmetry)}\nonumber \\
 & = & -R_{jikl}\,\,\,\,\,\,\,\,\,\,\,\,\,\text{(anti-symmetry in the first two indices)}\\
 & = & -R_{ijlk}\,\,\,\,\,\,\,\,\,\,\,\,\,\text{(anti-symmetry in the last two indices)}\nonumber
\end{eqnarray}

$\bullet$ The skew-symmetric property of the covariant Riemann tensor
with respect to the last two indices also applies to the mixed Riemann
tensor:
\begin{equation}
R_{\,\,jkl}^{i}=-R_{\,\,jlk}^{i}
\end{equation}

$\bullet$ As a consequence of the first and second anti-symmetric
properties of the covariant Riemann tensor, the entries of the Riemann
tensor with identical values of the first two indices or/and the last
two indices are zero.

$\bullet$ As a consequence of the skew-symmetric properties of the
Riemann tensor, all entries of the tensor with identical values of
more than two indices (e.g. $R_{iiji}$) are zero.

$\bullet$ In an $n$D space, the Riemann tensor has $n^{4}$ components.

$\bullet$ As a consequence of the symmetric and anti-symmetric properties
of the Riemann tensor, in an $n$D space there are three types of
distinct non-vanishing entries:

~~~~~~~A. Entries with only two distinct indices (type $R_{ijij}$)
which count:
\begin{equation}
N_{1}=\frac{n\left(n-1\right)}{2}
\end{equation}

~~~~~~~B. Entries with only three distinct indices (type $R_{ijik}$)
which count:
\begin{equation}
N_{2}=\frac{n\left(n-1\right)\left(n-2\right)}{2}
\end{equation}

~~~~~~~C. Entries with four distinct indices (type $R_{ijkl}$)
which count:
\begin{equation}
N_{3}=\frac{n\left(n-1\right)\left(n-2\right)\left(n-3\right)}{12}
\end{equation}

$\bullet$ By adding the numbers of the three types of non-zero distinct
entries, as given in the last point, it can be shown that the Riemann
tensor in an $n$D space has a total of
\begin{equation}
N_{\mathrm{RI}}=N_{1}+N_{2}+N_{3}=\frac{n^{2}\left(n^{2}-1\right)}{12}
\end{equation}
independent components which do not vanish identically. For example,
in a 2D Riemannian space the Riemann tensor has $2^{4}=16$ components;
however there is only one independent component (with the principal
suffix 1212) which is not identically zero represented by the following
four dependent components:
\begin{equation}
R_{1212}=R_{2121}=-R_{1221}=-R_{2112}
\end{equation}
Similarly, in a 3D Riemannian space the Riemann tensor has $3^{4}=81$
components but only six of these are distinct non-zero entries which
are the ones with the following principal suffixes:
\begin{equation}
1212,\,1313,\,1213,\,2123,\,3132,\,2323
\end{equation}
where the permutations of indices in each of these suffixes are subject
to the symmetric and anti-symmetric properties of the four indices
of the Riemann tensor, as in the case of a 2D space in the above example,
and hence these permutations do not produce independent entries.

$\bullet$ Following the pattern in the last point, in a 4D Riemannian
space the Riemann tensor has $4^{4}=256$ components but only 20 of
these are independent non-zero entries, while in a 5D Riemannian space
the Riemann tensor has $5^{4}=625$ components but only 50 are independent
non-zero entries.

$\bullet$ The Riemann tensor satisfies the following identity:
\begin{equation}
R_{ijkl;s}+R_{iljk;s}=R_{iksl;j}+R_{ikjs;l}
\end{equation}

$\bullet$ A necessary and sufficient condition that a manifold for
which there is a coordinate system with all the components of the
metric tensor being constants\footnote{This may be called flat pseudo Riemannian manifold.}
is that:
\begin{equation}
R_{ijkl}=0
\end{equation}

$\bullet$ On contracting the first covariant index with the contravariant
index of the Riemann tensor of the second kind we obtain:
\begin{equation}
\begin{aligned}R_{\,ikl}^{i} & =\partial_{k}\Gamma_{il}^{i}-\partial_{l}\Gamma_{ik}^{i}+\Gamma_{il}^{r}\Gamma_{rk}^{i}-\Gamma_{ik}^{r}\Gamma_{rl}^{i} & \,\,\,\,\,\, & \text{(\ensuremath{j=i} in Eq. \ref{eqRijklMix})}\\
 & =\partial_{k}\Gamma_{il}^{i}-\partial_{l}\Gamma_{ik}^{i}+\Gamma_{il}^{r}\Gamma_{rk}^{i}-\Gamma_{rk}^{i}\Gamma_{il}^{r} &  & \text{(relabeling dummy \ensuremath{i,r} in last term)}\\
 & =\partial_{k}\Gamma_{il}^{i}-\partial_{l}\Gamma_{ik}^{i}\\
 & =\partial_{k}\left[\partial_{l}\left(\ln\sqrt{g}\right)\right]-\partial_{l}\left[\partial_{k}\left(\ln\sqrt{g}\right)\right] &  & \text{(Eq. \ref{eqGammaaia})}\\
 & =\partial_{k}\partial_{l}\left(\ln\sqrt{g}\right)-\partial_{l}\partial_{k}\left(\ln\sqrt{g}\right)\\
 & =\partial_{k}\partial_{l}\left(\ln\sqrt{g}\right)-\partial_{k}\partial_{l}\left(\ln\sqrt{g}\right) &  & \text{(\ensuremath{C^{2}} condition is assumed)}\\
 & =0
\end{aligned}
\end{equation}
That is:
\begin{equation}
R_{\,ikl}^{i}=0
\end{equation}

\subsubsection{Bianchi identities}

$\bullet$ The Riemann tensor of the first and second kind satisfies
a number of identities called the Bianchi identities.

$\bullet$ The first Bianchi identity is:
\begin{equation}
\begin{aligned}R_{ijkl}+R_{iljk}+R_{iklj} & =0 & \,\,\,\,\,\,\,\,\,\,\,\,\,\,\,\,\,\,\,\,\,\, & \text{(first kind)}\\
R_{\,\,jkl}^{i}+R_{\,\,ljk}^{i}+R_{\,\,klj}^{i} & =0 &  & \text{(second kind)}
\end{aligned}
\end{equation}
These two forms of the first identity can be obtained from each other
by the raising and lowering operators.

$\bullet$ The above first Bianchi identity is an instance of the
fact that by fixing the position of one of the four indices and permuting
the other three indices cyclically, the algebraic sum of these three
permuting forms is zero, that is:
\begin{equation}
\begin{aligned}R_{ijkl}+R_{iljk}+R_{iklj} & =0 & \,\,\,\,\,\,\,\,\,\,\,\,\,\,\, & \text{(\ensuremath{i} fixed)}\\
R_{ijkl}+R_{ljik}+R_{kjli} & =0 &  & \text{(\ensuremath{j} fixed)}\\
R_{ijkl}+R_{likj}+R_{jlki} & =0 &  & \text{(\ensuremath{k} fixed)}\\
R_{ijkl}+R_{kijl}+R_{jkil} & =0 &  & \text{(\ensuremath{l} fixed)}
\end{aligned}
\end{equation}

$\bullet$ Another one of the Bianchi identities is:
\begin{equation}
\begin{aligned}R_{ijkl;m}+R_{ijlm;k}+R_{ijmk;l} & =0 & \,\,\,\,\,\,\,\,\,\,\,\,\,\,\, & \text{(first kind)}\\
R_{\,\,jkl;m}^{i}+R_{\,\,jlm;k}^{i}+R_{\,\,jmk;l}^{i} & =0 &  & \text{(second kind)}
\end{aligned}
\end{equation}
Again, these two forms can be obtained from each other by the raising
and lowering operators.

$\bullet$ The Bianchi identities are valid regardless of the metric.

\subsection{Ricci Tensor}

$\bullet$ The Ricci tensor of the first kind is obtained by contracting
the contravariant index with the last covariant index of the Riemann
tensor of the second kind, that is:
\begin{equation}
R_{ij}=R_{ija}^{a}=\partial_{j}\Gamma_{ia}^{a}-\partial_{a}\Gamma_{ij}^{a}+\Gamma_{bj}^{a}\Gamma_{ia}^{b}-\Gamma_{ba}^{a}\Gamma_{ij}^{b}
\end{equation}
and hence it is a rank-2 tensor.

$\bullet$ The Ricci tensor, as given by the last equation, can be
written in the following mnemonic determinantal form:
\begin{equation}
R_{ij}=\begin{vmatrix}\begin{array}{cc}
\partial_{j} & \partial_{a}\\
\Gamma_{ij}^{a} & \Gamma_{ia}^{a}
\end{array}\end{vmatrix}+\begin{vmatrix}\begin{array}{cc}
\Gamma_{bj}^{a} & \Gamma_{ba}^{a}\\
\Gamma_{ij}^{b} & \Gamma_{ia}^{b}
\end{array}\end{vmatrix}
\end{equation}

$\bullet$ Because of Eq. \ref{eqGammaaia} (i.e. $\Gamma_{ij}^{j}=\partial_{i}\left(\ln\sqrt{g}\right)$),
the Ricci tensor can also be written in the following forms as well
as several other forms:
\begin{eqnarray}
R_{ij} & = & \partial_{j}\partial_{i}\left(\ln\sqrt{g}\right)-\partial_{a}\Gamma_{ij}^{a}+\Gamma_{bj}^{a}\Gamma_{ia}^{b}-\Gamma_{ij}^{b}\partial_{b}\left(\ln\sqrt{g}\right)\nonumber \\
 & = & \partial_{j}\partial_{i}\left(\ln\sqrt{g}\right)+\Gamma_{bj}^{a}\Gamma_{ia}^{b}-\partial_{a}\Gamma_{ij}^{a}-\Gamma_{ij}^{b}\partial_{b}\left(\ln\sqrt{g}\right)\label{eqRij}\\
 & = & \partial_{j}\partial_{i}\left(\ln\sqrt{g}\right)+\Gamma_{bj}^{a}\Gamma_{ia}^{b}-\frac{1}{\sqrt{g}}\partial_{a}\left(\sqrt{g}\Gamma_{ij}^{a}\right)\nonumber
\end{eqnarray}
where $g$ is the determinant of the covariant metric tensor.

$\bullet$ The Ricci tensor of the first kind is symmetric, that is:
\begin{equation}
R_{ij}=R_{ji}
\end{equation}

$\bullet$ On raising the first index of the Ricci tensor of the first
kind, the Ricci tensor of the second kind is obtained:
\begin{equation}
R_{\,\,j}^{i}=g^{ik}R_{kj}
\end{equation}

$\bullet$ The Ricci scalar, which is also called the curvature scalar
and the curvature invariant, is the result of contracting the indices
of the Ricci tensor of the second kind, that is:
\begin{equation}
R=R_{\,i}^{i}
\end{equation}

$\bullet$ Since the Ricci scalar is obtained by raising a subscript
index of the Ricci tensor of the first kind using the raising operator
followed by contracting the two indices, it can be written as:
\begin{equation}
R=g^{ij}R_{ij}=g^{ij}\left[\partial_{j}\partial_{i}\left(\ln\sqrt{g}\right)+\Gamma_{bj}^{a}\Gamma_{ia}^{b}-\frac{1}{\sqrt{g}}\partial_{a}\left(\sqrt{g}\Gamma_{ij}^{a}\right)\right]
\end{equation}
where the expression in the square brackets is obtained from the last
line of Eq. \ref{eqRij}; similar expressions can be obtained from
the other lines of that equation.

$\bullet$ In an $n$D space, the Ricci tensor has $n^{2}$ entries.
However, because of its symmetry it possesses a maximum of
\begin{equation}
N_{\mathrm{RD}}=\frac{n\left(n+1\right)}{2}\label{eqRicciIndependentN}
\end{equation}
distinct entries. As an example, in the 4D manifold of general relativity
$n=4$, and hence the Ricci tensor has $4^{2}=16$ components. However,
due to the symmetry of the Ricci tensor there are only ten independent
entries according to Eq. \ref{eqRicciIndependentN}. The gravitational
field equations in a free space are obtained by setting the Ricci
tensor components equal to zero, and hence there are ten partial differential
equations describing the gravitational field in this space according
to the general relativistic mechanics.

\subsection{Einstein Tensor}

$\bullet$ The Einstein tensor $\mathbf{G}$ is a rank-2 tensor defined
in terms of the Ricci tensor $\mathbf{R}$ and the Ricci curvature
scalar $R$ as follow:\footnote{We notate this tensor with $\mathbf{G}$ rather than $\mathbf{E}$,
which is more natural, to avoid potential confusion with the first
displacement gradient tensor (see $\S$ \ref{subDisplacementGradientTensors}).}
\begin{equation}
\begin{aligned}G_{mn} & =R_{mn}-\frac{1}{2}g_{mn}R & \,\,\,\,\,\,\,\,\,\,\,\,\,\,\,\,\, & \text{(covariant)}\\
G^{mn} & =R^{mn}-\frac{1}{2}g^{mn}R &  & \text{(contravariant)}\\
G_{n}^{m} & =R_{n}^{m}-\frac{1}{2}\delta_{n}^{m}R &  & \text{(mixed)}
\end{aligned}
\end{equation}

$\bullet$ Since both the Ricci tensor and the metric tensor are symmetric,
the Einstein tensor is symmetric as well.

$\bullet$ The divergence of the Einstein tensor vanishes at all points
of the space for any Riemannian metric.

$\bullet$ On contracting the Bianchi identity twice with using the
anti-symmetric properties of the Riemann tensor we obtain:
\begin{equation}
G_{\,\,\,;n}^{mn}=0
\end{equation}
which is inline with the above statement. The following form can also
be derived based on the Bianchi identity:
\begin{equation}
g^{jk}G_{ki;j}=0
\end{equation}

\subsection{Infinitesimal Strain Tensor\label{subInfinitesimalStrainTensor}}

$\bullet$ This is a rank-2 tensor which describes the state of strain
in a continuum medium and hence it is used in continuum and fluid
mechanics.

$\bullet$ The infinitesimal strain tensor $\boldsymbol{\gamma}$
is defined by:\footnote{Some authors do not include the factor $\frac{1}{2}$ in the definition
of $\boldsymbol{\gamma}$.}
\begin{equation}
\boldsymbol{\gamma}=\frac{\nabla\mathbf{d}+\nabla\mathbf{d}^{T}}{2}
\end{equation}
where $\mathbf{d}$ is the displacement vector and the superscript
$T$ represents matrix transposition. The displacement vector $\mathbf{d}$
represents the change in distance and direction which an infinitesimal
element of the medium experiences as a consequence of the applied
stress.

$\bullet$ In Cartesian coordinates with tensor notation, the last
equation is given as:
\begin{equation}
\gamma_{ij}=\frac{\partial_{i}d_{j}+\partial_{j}d_{i}}{2}
\end{equation}

\subsection{Stress Tensor\label{subStressTensor}}

$\bullet$ The stress tensor, which is also called Cauchy stress tensor,
is a rank-2 symmetric\footnote{In fact it is symmetric in many applications (e.g. in the flow of
Newtonian fluids) but not all, as it can be asymmetric in some cases.
We also choose to define it within the context of Cauchy stress law
which is more relevant to the continuum mechanics; however it can
be defined differently in other disciplines and in a more general
form.} tensor used for \textit{transforming} a normal vector to a surface
\textit{to} a traction vector acting on that surface, that is:
\begin{equation}
\mathbf{T}=\boldsymbol{\sigma}\mathbf{n}
\end{equation}
where $\mathbf{T}$ is the traction vector, $\boldsymbol{\sigma}$
is the stress tensor and $\mathbf{n}$ is the normal vector. This
is usually represented in tensor notation using Cartesian coordinates
as:
\begin{equation}
T_{i}=\sigma_{ij}n_{j}
\end{equation}

$\bullet$ The diagonal components of the stress tensor represent
normal stresses while the off-diagonal components represent shear
stresses.

$\bullet$ Because the stress tensor is symmetric, in an $n$D space
it possesses $\frac{n\left(n+1\right)}{2}$ independent components
instead of $n^{2}$. Hence in a 3D space (which is the ordinary space
for this tensor) it has six independent components.

$\bullet$ In fluid dynamics, the stress tensor (or total stress tensor)
is decomposed into two main parts: a viscous contribution part and
a pressure contribution part. The viscous part may then be split into
a normal stress and a shear stress while the pressure part may be
split into a hydrostatic pressure and an extra pressure.

\subsection{Displacement Gradient Tensors\label{subDisplacementGradientTensors}}

$\bullet$ These are rank-2 tensors which are denoted by $\mathbf{E}$
and $\boldsymbol{\Delta}$. They are defined in Cartesian coordinates
using tensor notation as:
\begin{equation}
E_{ij}=\frac{\partial x_{i}}{\partial x_{j}^{'}}\,\,\,\,\,\,\,\,\,\,\,\,\,\,\,\,\,\,\,\,\,\,\,\&\,\,\,\,\,\,\,\,\,\,\,\,\,\,\,\,\,\,\,\,\,\,\,\Delta_{ij}=\frac{\partial x_{i}^{'}}{\partial x_{j}}
\end{equation}
where $x$ and $x^{'}$ represent the Cartesian coordinates of an
observed continuum particle at the present and past times respectively.
These tensors may also be called deformation gradient tensors.

$\bullet$ $\mathbf{E}$ quantifies the displacement of the particle
at the present time relative to its position at the past time, while
$\boldsymbol{\Delta}$ quantifies its displacement at the past time
relative to its position at the present time.

$\bullet$ From their definitions, it is obvious that $\mathbf{E}$
and $\boldsymbol{\Delta}$ are inverses of each other and hence:
\begin{equation}
E_{ik}\Delta_{kj}=\delta_{ij}
\end{equation}

\subsection{Finger Strain Tensor}

$\bullet$ This, which may also be called the left Cauchy-Green deformation
tensor, is a rank-2 tensor used in the fluid and continuum mechanics
to describe the strain in a continuum object, e.g. fluid, in a series
of time frames. It is defined as:
\begin{equation}
\mathbf{B}=\mathbf{E}\cdot\mathbf{E}^{T}
\end{equation}
which in Cartesian coordinates with tensor notation becomes:
\begin{equation}
B_{ij}=\frac{\partial x_{i}}{\partial x_{k}^{'}}\frac{\partial x_{j}}{\partial x_{k}^{'}}
\end{equation}
where $\mathbf{E}$ is the first displacement gradient tensor as defined
in $\S$ \ref{subDisplacementGradientTensors}, the superscript $T$
represents matrix transposition, and the indexed $x$ and $x^{'}$
represent the Cartesian coordinates of an element of the continuum
at the present and past times respectively.

\subsection{Cauchy Strain Tensor}

$\bullet$ This, which may also be called the right Cauchy-Green deformation
tensor, is the inverse of the Finger strain tensor and hence it is
denoted by $\mathbf{B}^{-1}$. Consequently, it is defined as:
\begin{equation}
\mathbf{B}^{-1}=\boldsymbol{\Delta}^{T}\cdot\boldsymbol{\Delta}
\end{equation}
which in Cartesian coordinates with tensor notation becomes:
\begin{equation}
B_{ij}^{-1}=\frac{\partial x_{k}^{'}}{\partial x_{i}}\frac{\partial x_{k}^{'}}{\partial x_{j}}
\end{equation}
where $\boldsymbol{\Delta}$ is the second displacement gradient tensor
as defined in $\S$ \ref{subDisplacementGradientTensors}.

$\bullet$ The Finger and Cauchy strain tensors may be labeled as
``finite strain tensors'' as opposite to infinitesimal strain tensors.
They are symmetric positive definite tensors; moreover they become
the unity tensor when the change in the state of the object from the
past to the present times consists of rotation and translation with
no deformation.

\subsection{Velocity Gradient Tensor}

$\bullet$ This is a rank-2 tensor which is often used in fluid dynamics
and rheology. As its name suggests, it is the gradient of the velocity
vector $\mathbf{v}$ and hence it is given in Cartesian coordinates
by:
\begin{equation}
\left[\nabla\mathbf{v}\right]_{ij}=\partial_{i}v_{j}
\end{equation}

$\bullet$ The velocity gradient tensor in other coordinate systems
can be obtained from the expressions of the gradient of vectors in
these systems, as given, for instance, in $\S$ \ref{subCylindricalSystem}
and $\S$ \ref{subSphericalSystem} for cylindrical and spherical
coordinates.

$\bullet$ The term ``velocity gradient tensor'' my also be used
for the transpose of this tensor, i.e. $\left(\nabla\mathbf{v}\right)^{T}$.

$\bullet$ The velocity gradient tensor is usually decomposed into
a symmetric part which is the rate of strain tensor $\mathbf{S}$
(see $\S$ \ref{subRateofStrain}), and an anti-symmetric part which
is the vorticity tensor $\mathbf{\bar{S}}$ (see $\S$ \ref{subVorticityTensor}),
that is:
\begin{equation}
\nabla\mathbf{v}=\mathbf{S}+\mathbf{\bar{S}}
\end{equation}

\subsection{Rate of Strain Tensor\label{subRateofStrain}}

$\bullet$ This tensor, which is also called the rate of deformation
tensor, is the symmetric part of the velocity gradient tensor and
hence is given by:\footnote{Some authors do not include the factor $\frac{1}{2}$ in the definition
of $\mathbf{S}$ and $\mathbf{\bar{S}}$ and hence this factor is
moved to the definition of $\nabla\mathbf{v}$. Also these tensors
are commonly denoted by $\boldsymbol{\dot{\gamma}}$ and $\boldsymbol{\omega}$
respectively.}
\begin{equation}
\mathbf{S}=\frac{\nabla\mathbf{v}+\left(\nabla\mathbf{v}\right)^{T}}{2}
\end{equation}
which, in tensor notation with Cartesian coordinates, is given by:
\begin{equation}
S_{ij}=\frac{\partial_{i}v_{j}+\partial_{j}v_{i}}{2}
\end{equation}

$\bullet$ The rate of strain tensor is a quantitative measure of
the local rate at which neighboring material elements of a deforming
continuum move with respect to each other.

$\bullet$ As a rank-2 symmetric tensor, it has $\frac{n\left(n+1\right)}{2}$
independent components which is six in a 3D space.

$\bullet$ The rate of strain tensor is related to the infinitesimal
strain tensor (refer to $\S$ \ref{subInfinitesimalStrainTensor})
by:
\begin{equation}
\mathbf{S}=\frac{\partial\boldsymbol{\gamma}}{\partial t}
\end{equation}
where $t$ is time. Hence, the rate of strain tensor is normally denoted
by $\boldsymbol{\dot{\gamma}}$ where the dot represents the temporal
rate of change.

\subsection{Vorticity Tensor\label{subVorticityTensor}}

$\bullet$ This is the anti-symmetric part of the velocity gradient
tensor and hence is given by:
\begin{equation}
\mathbf{\bar{S}}=\frac{\nabla\mathbf{v}-\left(\nabla\mathbf{v}\right)^{T}}{2}
\end{equation}
which, in tensor notation with Cartesian coordinates, is given by:
\begin{equation}
\bar{S}_{ij}=\frac{\partial_{i}v_{j}-\partial_{j}v_{i}}{2}
\end{equation}

$\bullet$ The vorticity tensor quantifies the local rate of rotation
of a deforming continuum medium.

$\bullet$ As a rank-2 anti-symmetric tensor, it has $\frac{n\left(n-1\right)}{2}$
independent components which is three in a 3D space. These three components
added to the six components of the rate of strain tensor give nine
independent components which is the total number of independent components
of their parent tensor $\nabla\mathbf{v}$.

\pagebreak{}

\phantomsection
\addcontentsline{toc}{section}{References}
\bibliographystyle{unsrt}

\end{document}